\documentstyle[amssymb,12pt]{amsart}%\input amssym.def
\newtheorem{prop}{Proposition}[section]
\newtheorem{prop:def}{Proposition-Definition}[section]

\newtheorem{lemma}{Lemma}[section]
\newtheorem{thm}{Theorem}[section]
\newtheorem{cor}{Corollary}[section]

\theoremstyle{remark}
\newtheorem{remark}{Remark}

\begin{document}
\newcommand{\nc}{\newcommand} \nc{\on}{\operatorname}
\nc{\pa}{\partial} \nc{\cA}{{\cal A}}\nc{\cB}{{\cal B}}\nc{\cC}{{\cal
    C}} \nc{\cE}{{\cal E}}\nc{\cG}{{\cal G}}\nc{\cH}{{\cal H}}
\nc{\cX}{{\cal X}}\nc{\cR}{{\cal R}}\nc{\cL}{{\cal L}} \nc{\cK}{{\cal
    K}} \nc{\sh}{\on{sh}}\nc{\val}{\on{val}}
\nc{\Id}{\on{Id}}\nc{\Diff}{\on{Diff}}
\nc{\ad}{\on{ad}}\nc{\Der}{\on{Der}}\nc{\End}{\on{End}}
\nc{\res}{\on{res}}\nc{\ddiv}{\on{div}}
\nc{\card}{\on{card}}\nc{\dimm}{\on{dim}}
\nc{\Jac}{\on{Jac}}\nc{\Ker}{\on{Ker}}
\nc{\Imm}{\on{Im}}\nc{\limm}{\on{lim}}\nc{\Ad}{\on{Ad}}
\nc{\ev}{\on{ev}} \nc{\Hol}{\on{Hol}}\nc{\Det}{\on{Det}}
\nc{\Bun}{\on{Bun}}\nc{\diag}{\on{diag}}
\nc{\de}{\delta}\nc{\si}{\sigma}\nc{\ve}{\varepsilon}
\nc{\al}{\alpha}\nc{\vp}{\varphi} \nc{\CC}{{\mathbb
    C}}\nc{\ZZ}{{\mathbb Z}} \nc{\VV}{{\mathbb V}} \nc{\nn}{{\mathbf
    n}} \nc{\xxi}{{\tilde\xi}} \nc{\NN}{{\mathbb N}}\nc{\zz}{{\mathbf
    z}} \nc{\AAA}{{\mathbb A}}\nc{\cO}{{\cal O}} \nc{\cF}{{\cal
    F}}\nc{\cM}{{\cal M}} \nc{\la}{{\lambda}}\nc{\G}{{\mathfrak
    g}}\nc{\mm}{{\mathfrak m}} \nc{\A}{{\mathfrak a}}
\nc{\HH}{{\mathfrak h}} \nc{\N}{{\mathfrak n}}\nc{\B}{{\mathfrak b}}
\nc{\La}{\Lambda}
\nc{\g}{\gamma}\nc{\eps}{\epsilon}\nc{\wt}{\widetilde}
\nc{\wh}{\widehat} \nc{\bn}{\begin{equation}}\nc{\en}{\end{equation}}
\nc{\SL}{{\mathfrak{sl}}}\nc{\ttt}{{\mathfrak{t}}}

% ****** GISPIC **********
%
%** by GISLI MASON *******
%
%**for commutative diagrams
%

\newcommand{\ldar}[1]{\begin{picture}(10,50)(-5,-25)
\put(0,25){\vector(0,-1){50}}
\put(5,0){\mbox{$#1$}} 
\end{picture}}

\newcommand{\lrar}[1]{\begin{picture}(50,10)(-25,-5)
\put(-25,0){\vector(1,0){50}}
\put(0,5){\makebox(0,0)[b]{\mbox{$#1$}}}
\end{picture}}

\newcommand{\luar}[1]{\begin{picture}(10,50)(-5,-25)
\put(0,-25){\vector(0,1){50}}
\put(5,0){\mbox{$#1$}}
\end{picture}}

\title[Commuting difference operators associated to complex curves]
{Commuting differential and difference operators associated to complex
  curves, II}

\author{B. Enriquez}

\address{B.E.: Centre de Math\'ematiques, Ecole Polytechnique, 
UMR 7640 du CNRS, 91128 Palaiseau, France}

\address{FIM, ETH-Zentrum, HG G46, CH-8092 Zurich, Switzerland}

\author{G. Felder}

\address{G.F.: D-Math, ETH-Zentrum, HG G44, CH-8092 Zurich,
  Switzerland}

\date{December 1998}

% \begin{abstract}
% \end{abstract}

\maketitle

\subsection*{Introduction}

This paper is a sequel to \cite{commI}. Our main aim is to construct a
commuting family of difference-evaluation operators $(T_z^{(\Pi)})_z$,
deforming the difference-evaluation operators $T_z^{class}$ of
\cite{commI}, and to interpret them as the action of the center of a
quantum algebra in the space of intertwiners of a ``regular''
subalgebra.

Let us recall first some points of \cite{commI}. In that paper, we
proposed a functional approach to the Knizhnik-Zamolodchikov-Bernard
(KZB) connection, relying on the functional picture for conformal
blocks of \cite{FS}. Recall that conformal blocks are associated to a
complex curve $X$ with a marked point $P_0$, a simple Lie algebra
$\bar\G$ and representations $\VV$ and $V$ of $\G$ and $\G^{out}$,
where $\G$ is the Kac-Moody algebra $(\bar\G \otimes \cK)\oplus \CC
K$, and $\G^{out}$ is the Lie subalgebra of $\G$ formed of the
currents regular outside $P_0$ (we denote by $\cK$ is the local field
of $X$ at $P_0$), and defined as the space of $\G^{out}$-intertwiners
$\psi$ from $\VV$ to $V$. Twisted conformal blocks are defined in the
same way, replacing $\G^{out}$ by the Lie subalgebra
$\G^{out}_{\la_0}$ of $\G$ formed of the maps $x$ from the universal
cover of $X$, regular outside the preimage of $P_0$, with
transformation properties $x(\gamma_{A_a}z) = x(z)$ and
$x(\gamma_{B_a}z) = e^{\la^{(0)}_a}x(z) e^{-\la^{(0)}_a}$, where
$\gamma_{A_a}$ and $\gamma_{B_a}$ are deck transformations
corresponding to $a$- and $b$-cycles; $\la_0 = (\la^{(0)}_a)$ belongs
to $\bar\HH^g$, where $\bar\HH$ is the Cartan subalgebra of $\bar\G$
and $\G^{out}_{\la_0}$ is the Lie subalgebra of $\G$.  We parametrize
the space of twisted conformal blocks by associating to $\psi_{\la_0}$
the twisted correlation functions of currents of $\G$ associated to
the simple root generators of its nilpotent subalgebra $\bar\N_+$.
Denote by $e_i,f_i,h_i$ the Chevalley generators of $\bar\G$, so that
the $e_i$ generate $\bar\N_+$, and set $x[f] = x \otimes f$, for $x$
in $\bar\G$, $f$ in $\cK$. To a vector $v$ of $\VV$, annihilated by
the $h_i[z^{k'}],f_i[z^{1-g +k}]$, $k'>0, k\geq 0$ (this property is
shared by the extremal vectors in integrable modules), and to an
intertwiner $\psi_{\la_0}$, weon $\VV$, $\G^{out}_{\la_0}$-invariant,
we associate the generating series
$$ f(\la^{(i)}_a | z^{(i)}_j) = \langle \psi_{\la} [ \prod_i
\prod_{j=1}^{n_j} e_i(z^{(i)}_j) v ], \xi \rangle ,
$$ where $\psi_\la = \psi_{\la_0} \circ e^{\sum_a (\la_a -
  \la_a^{(0)}) h[r_a]}$, the $r_a$ are multivalued functions on $X$,
constant along $a$-cycles and with additive constants along
$b$-cycles, and $\xi$ is a lowest weight form on $V$. The $\la_a$ are
formal parameters near $(\la_a^{(0)})$ and the $z_i$ are formal
parameters near $P_0$.

In \cite{commI}, we expressed the KZB connection in terms of these
correlation functions. Let $T(z)$ denote the Sugawara tensor; it is a
series in $(U\G)_{loc}[[z,z^{-1}]]$, where $(U\G)_{loc}$ is the local
completion of the universal enveloping algebra $U\G$. In the case
$\bar\G = \SL_2$, we have

\begin{thm} (see \cite{commI})
  Let $T_z^{class}$ be the differential-evaluation operator acting on
  functions $f(\la_1,...,\la_g| z_1, ..., z_n)$ as
\begin{align*}
  & T_z^{class} = {1\over 2} [\sum_a \omega_a(z) \pa_{\la_a} + 2
  \sum_{i} G^{(I)} (z,z_i) - \sum_{j} \La_j G^{(I)}(z,P_j) ]^2 \\ & +
  \sum_a D_z^{(2\la)}\omega_a(z) \pa_{\la_a} + 2\sum_i D_z^{(2\la)}
  G^{(I)}(z,z_i) - \sum_j \La_j D_z^{(2\la)}
  G^{(I)}(z,P_j) + k \omega_{2\la}(z) \\ & + \sum_{i=1}^n \left( -2
    G^{(I)}_{2\la}(z,z_i)[\sum_a \omega_a(z_i) \pa_{\la_a} + 2
    \sum_{j\neq i} G^{(I)}(z_i,z_j)- \sum_{k} \La_k G^{(I)}(z_i,P_k)]
  \right.  \\ & \left. -4 G^{(I)}_{2\la}(z,z_i)G^{(I)}(z_i,z) + 2 k
    d_{z_i}G^{(I)}_{2\la}(z,z_i) \right) \circ \ev_z^{(i)},
\end{align*} 
where 
$$ (\ev_z^{(i)}f)(\la| z_1, ... , z_n) = f(\la|z_1, ..., z, ...,
z_n)
$$ ($z$ in $i$th position), we set $\la = (\la_1,...,\la_g)$,
$D_z^{(2\la)}$ is a connection on the bundle $K$ of differentials on
$X$ and has simple pole at $P_0$, the $\omega_a$ are the holomorphic
one-forms associated with the $a$-cycles, $\omega_{2\la}$ is a
quadratic differential with double poles at $P_0$, $G^{(I)}$ and
$G^{(I)}_\la$ are (twisted) Green functions, $P_j$ are some point of
$X - \{P_0\}$ and $\La_i$ are some numbers.

If $V$ is the product $\otimes_i V_{-\La_i}(P_i)$ of evaluation
modules ($V_{-\La_i}$ is the $\SL_2$-module with lowest weight
$-\La_i$), we have the equality
$$ \langle \psi_\la [T(z) \prod_i e_{i}(z^{(i)}_j)v], \xi \rangle =
(T_z^{class}f)(\la|z_1,...,z_n),
$$ if $k$ is the level of $\VV$.  The operators $T_z^{class}$ commute
when $k = -2$.
\end{thm}
When $X$ is $\CC P^1$, the expression of $T_z^{class}$ is similar to
the expression for the action of the Hamiltonians on Bethe vectors
obtained in the Bethe ansatz approach to the Gaudin system (see
\cite{FFR}).

\medskip 

In the present paper, we repeat these steps of \cite{commI} in the
quantum case, at the critical level. We replace the Kac-Moody algebra
$\G$ by the quantum group $U_{\hbar,\omega}\G$ associated to a pair
$(X,\omega)$ of a curve $X$ and a rational differential $\omega$
(\cite{Ann3}). The relations for this algebra depend on the choice of
a Lagrangian subspace of $\cK$, that we construct in sect.\ 
\ref{sect:geom}. We recall the presentation of $U_{\hbar,\omega}\G$ in
terms of generating fields $e(z),f(z),k^{\pm}(z)$ (sect.\ 
\ref{pres:G}). The algebra $U_{\hbar,\omega}\G$ contains a subalgebra
$U_\hbar\G^{out}$, which is a flat deformation of the enveloping
algebra of $\G^{out}$ (\cite{ER:qH}).

Let $\bar\G = \bar\N_+ \oplus \bar\HH \oplus \bar\N_-$ be the Cartan
decomposition of $\bar\G$. Let $\mm$ be the maximal ideal at $P_0$ and
$\B_{in}$ be the subalegbra of $\G$ defined as $\B_{in} = (\bar\HH
\otimes \mm) \oplus (\bar\N_+ \otimes \cK)$.  We construct, in
$U_{\hbar,\omega}\G$, a subalgebra isomorphic to $(U\B_{in})[[\hbar]]$
(sect.\ \ref{sub+}).  This subalgebra is expressed in terms of ``new''
generating fields $\wt e(z),\wt f(z)$ and $k_{tot}^{\pm}(z)$; we study
their relations in sect.\ \ref{other:pres}.  In the rational case,
such generating fields appeared in \cite{Sevost}.  We express a
generating function for central elements $T(z)$ deforming the Sugawara
tensor by the formula (see Thm.\ \ref{Sugawara})
$$ T(z) = : e(z) \wt f(z) :_{\la}+ a_\la(z) k^{+}_{tot}(z) + b_\la(z)
k^{-}_{tot}(z) ,
$$ where $:a(z)b(z):_\la$ denotes a normal ordered product, depending
on $\la$ and $a_\la(z)$ and $b_\la(z)$ are formal series of
$\cK[[\la_a - \la_a^{(0)}]][[\hbar]]$ defined by (\ref{a:la}) and
(\ref{b:la}). We also obtain another expression for $T(z)$ of the type
obtained in \cite{FR,RS}, see (\ref{old:T}).

We construct a subalgebra $U_\hbar \G^{out}_{\la_0}$ of
$U_{\hbar,\omega}\G$ in sect.\ \ref{subalg} deforming the enveloping
algebra of $\G^{out}_{\la_0}$ and study a class of its representations
(sect.\ \ref{fd:reps}). We show that such representations have a
lowest weight form $\xi$, such that
$$
\xi \circ f[r_{-2\la_0}] = 0, \quad \xi \circ k^+(z) = \pi(z) \xi, 
$$ for $r_{-2\la_0}$ in $R_{-2\la_0}$ and $\pi(z)$ a formal series,
which is an analogue of the Drinfeld polynomial.

To a module $\VV$ over $U_{\hbar,\omega}\G$, and to a morphism
$\psi_{\la_0}: \VV \to V$ of $U_{\hbar}\G^{out}_{\la_0}$-modules,
where $V$ is a product of evaluation modules, we associate the
correlation function
$$ f(\la_1,...,\la_g| z_1,...,z_n) = \langle \psi_{\la}[\wt e(u_1)
\cdots \wt e(u_n) v], \xi \rangle,
$$ where $\xi$ is a lowest weight form on $V$ and $\psi_\la =
\psi_{\la_0} \circ e^{\sum_i (\la_a - \la_a^{(0)})h[r_a]}$. We study
the functional properties of $f(\la_1,...,\la_g | z_1,...,z_g)$ in
sect.\ \ref{tw:corr}.

Our main result is then

\begin{thm} (see Thm.\ \ref{comm:ops})
  Let for any formal series $\Pi$, $(T_z^{(\Pi)})_z$ be the family of
  operators acting on functions $f(\la| u_1, ... , u_n)$, defined as
\begin{align*} 
  & T_z^{(\Pi)} = \Pi(z) a'_\la(z| u_1,...,u_n) \circ e^{\sum_a
    \omega'_a(z) \pa / \pa\la_a} + \Pi(q^{-\pa}z)^{-1} a''_\la(z|
  u_1,...,u_n) \circ e^{\sum_a \omega''_a(z) \pa / \pa\la_a} \\ & +
  \sum_i \Pi(u_i) c_\la^{\prime (i)} (z| u_1,...,u_n) \circ e^{\sum_a
    \omega'_a(z) \pa / \pa\la_a} \circ \ev_z^{(i)} \\ & + \sum_i
  \Pi(q^{-\pa}u_i)^{-1} c_\la^{\prime\prime (i)} (z| u_1,...,u_n)
  \circ e^{\sum_a \omega''_a(z) \pa / \pa\la_a} \circ \ev_z^{(i)},
\end{align*}
where the multiplication operators are denoted as functions, and we
set
$$ a'_\la(z| u_1,...,u_n) = a_\la(z)\prod_i q_\mm(z,u_i), \ 
a''_\la(z| u_1,..,u_n) = b'_\la(z)\kappa(z)
\prod_i q_\mm(q^{-\pa}z,u_i)^{-1} ,
$$
$$ c^{\prime (i)}_\la(z| u_1,...,u_n) = -{1\over\hbar}
G_{2\la}(z,u_i)q_\mm(u_i,z)\prod_{j\neq i} q_\mm(u_i,u_j),$$
$$ c^{\prime \prime(i)}_\la(z| u_1,..,u_n) = {1\over\hbar}
G_{2\la}(z,q^{-\pa}u_i)
\kappa(u_i)q_\mm(q^{-\pa}u_i,z)^{-1}\prod_{j\neq i}
q_\mm(q^{-\pa}u_i,u_j)^{-1},
$$
$$ \omega'_a) = \hbar {1\over{1+q^{-\pa}}}(\omega_a / \omega)(z) ,
\quad \omega''_a) = - \hbar {1\over{1+q^{\pa}}}(\omega_a / \omega)(z),
\quad G_{2\la}(z,w)= G^{(I)}_{2\la}(z,w)/\omega(z),
$$ where $\pa$ is the derivation associated with $\omega$, so that
$\pa f = df/\omega(z)$, $a_\la$, $b'_\la$, $q_\mm$ and $\kappa$ are
defined in (\ref{a:la}), (\ref{b':la}), (\ref{q:mm}) and
(\ref{kappa}), and $q = e^\hbar$. The operators $T_z^{(\Pi)}$ commute
and normalize first order difference operators $\hat f[\rho]$ defined
by (\ref{f:rho}).  Moreover, we have, if the subalgebra
$U_\hbar\B^{\geq 1-g}$ of $U_\hbar\B_{in}$ acts on $v$ by the
character $\chi_n$ (see sect.\ \ref{v;chi:n}),
$$ \langle \psi_\la[T(z) \wt e(u_1) \cdots \wt e(u_n) v ] , \xi
\rangle = T_z\{ \langle \psi_\la[\wt e(u_1) \cdots \wt e(u_n) v ] ,
\xi \rangle \}, 
$$ where $\Pi$ can be expressed in terms of $\pi$. We also set
$\Pi(q^\pa z) = (q^\pa\Pi)(z)$.
\end{thm}

The $T_z^{(\Pi)}$ are difference deformations of the $T_z^{class}$.
In the rational case, we identify the operators $T_z^{(\Pi)}$ with the
commuting family of operators provided by the Yangian action on the
hypergeometric spaces of \cite{Tar-Var} (see sect.\ \ref{hypergeom}).
In the elliptic case, we identify $T_z^{(\Pi = 1)}$ with the first
$q$-Lam\'e operator (rem.\ \ref{elliptic}).

Let us say some words about possible prolongations of the present
work:

1) noncritical level.  One could try to prove analogues of the
theta-behavior results of \cite{commI} for the twisted correlation
functions of integrable modules over $U_{\hbar,\omega}\G$. Another
problem is to find analogues of the KZB flows for noncritical level,
by extending the approach of \cite{KS} to the quasi-Hopf situation.

%b) $q$-deformation of the $\pa_{z_i}$ flows out of critical level. One
%may try to use for this the Sugawara action approach: in case the
%$i$th module is a dim $2$ evaluation module, evaluation at $z = z_i$
%of the difference operator $T(z)$ gives exactly the $i$th qKZ
%operator. One may also try to use a construction \`a la \cite{KS}. At
%zero level, this would give a commuting family of operators on the
%$(\otimes_i V_i)^B$ ....  Understanding the elliptic qKZB in this way.

2) versions where $\hbar$ takes complex values. When $\omega$ is the
pull-back of the form $dz$ or $dz/z$ from a morphism $X\to \CC P^1$ or
$X \to E$, $E$ some elliptic curve, $q^\pa$ at least makes sense as
some correspondence on $X$. It could then be possible to find a
presentation of $U_{\hbar,\omega}\G$ allowing for complex values of
$\hbar$. This was done in \cite{examples} in the case $X = \CC P^1$,
$\omega = z^N dz$.

%$q^\pa$ would then be replaced by some set of automorphisms of $X$ and
%$q(z,w)$ by some rational function on $X^2$...

3) Bethe ansatz for the operators $T_z$. In \cite{FFR}, Bethe
equations for the Gaudin system are shown to be equivalent to the
existence of intertwining operators at critical level, and in turn to
a trivial monodromy condition for some connection. The similar study
should be possible for the systems constructed here, so that they
could be viewed as $q$-deformations of the differential systems
arising in \cite{BD}.

\medskip 

We would like to express our thanks to B.\ Feigin, E.\ Frenkel and V.\ 
Tarasov for discussions about this paper; the first author would like
to express his gratitude to A.-S.\ Sznitman for invitation to the FIM,
ETHZ, where this work was done.

\section{Geometric setting} \label{sect:geom}

\subsection{Isotropic supplementaries} 
Let $X$ be a smooth compact complex curve of genus $g$, endowed with a
nonzero holomorphic form $\omega$. Let $\sum_{i=1}^p n_i P_i$ be the
divisor of $\omega$ (we have $n_i >0$, $\sum_i n_i = 2(g-1)$). Let for
each $i$, $\cK_i$ be the local field at $P_i$, $\cO_i$ the local ring
at this point and $\mm_i$ the maximal ideal of $\cO_i$. Define $\cK$
as $\oplus_{i=1}^p \cK_i$ and $R$ as the space of rational functions
on $X$, regular outside $\{P_i\}$; we view it as a subring of $\cK$.
For each $i$, let $z_i$ be a local coordinate at $P_i$. Then $\cO_i =
\CC[[z_i]]$, $\mm_i = z_i \cO_i$ and $\cK_i =\CC((z_i))$.

$\cK$ is endowed with a scalar product $\langle , \rangle_\cK$ defined
by
$$ \langle f,g \rangle_\cK = \sum_{i=1}^p \res_{P_i}(fg\omega).
$$ Let us fix on $X$ a choice of $a$- and $b$-cycles $(A_a)_{1\leq a
  \leq g}$ and $(B_a)_{1\leq a \leq g}$. Let $\wt X$ be the universal
cover of $X$ and $\pi: \wt X \to X$ be the cover map. Denote by
$\gamma_{A_a}$ and $\gamma_{B_a}$ the deck transformations associated
with the cycles $A_a$ and $B_a$. 

\begin{lemma}
  There exists a linearly independent family of $R$ formed by
% the $\omega_i /\omega$ and of 
  elements $f_{(m_i)}$, where $m_i, i = 1, \ldots,p$ are integers such
  that $m_i \geq n_i$ for each $i$ and $\sum_i m_i \geq \sum_i n_i
  +2$, with $\on{val}_{P_i}(f_{(m_i)}) = - m_i$, for each $i = 1,
  \ldots,p$.
\end{lemma}

{\em Proof.} Let us first construct the $f_{(m_i)}$. Assume $m_j \geq
n_j + 1$, then by the Riemann-Roch theorem,
$$ h^0(\cO(\sum_i m_i P_i)) - h^0(\cO(\sum_i m_i P_i - P_j)) = 1 +
h^1(\cO(\sum_i m_i P_i)) - h^1(\cO(\sum_i m_i P_i - P_j)) ; 
$$
by Serre duality this is equal to 
$$ 1 + h^0(\cO(\sum_i (n_i - m_i) P_i)) - h^0(\cO(\sum_i (n_i - m_i)
P_i + P_j)). 
$$ All the $n_i - m_i + \delta_{ij}$ are $\leq 0$ and their sum is
$<0$, so not all of them are zero. Therefore both $h^0$ vanish. This
proves the existence of the $f_{(m_i)}$.  \hfill \qed \medskip

\begin{lemma}
  We have $g$ functions $r_a$ defined on $\wt X$, regular outside
  $\pi^{-1}(\{P_i\})$, such that

i) $r_a \circ \gamma_{B_b} = r_a -
  \delta_{ab}$, 

ii) $\on{val}_{P_i}(r_a) \geq -n_i - \delta_{i1}$ and

iii) $\int_{A_a} r_b \omega = {1\over 2} \int_{A_a} \omega \delta_{ab}$.   
\end{lemma}

{\em Proof.} The existence of rational functions $\tilde r_a$ defined
on $\wt X$, regular outside $\pi^{-1}(\{P_i\})$ and satisfying {\em
  i)} is a consequence of \cite{commI} Cor.\ 1.1. Adding to them
suitable combinations of the $f_{(m_i)}$, one gets functions $\bar
r_a$ satisfying both {\em i)} and {\em ii)}. Let $\omega_a$ be a basis
of the space of holomorphic one-froms on $X$. The ratios
$\omega_a/\omega$ are elements of $R$, with valuation at each $P_i$
less or equal to $-n_i$. Adding to the $\bar r_a$ suitable
combinations of the $\omega_a / \omega$, one obtains elements $r_a$
satisfying {\em i)}, {\em ii)} and {\em iii)}.  \hfill \qed \medskip

\begin{prop} \label{Lambda}
  Set $\Lambda = (\oplus_a \CC r_a) \oplus (\mm_1 \oplus \cO_2 \oplus
  \cdots \oplus \cO_p)$. We have a direct sum decomposition $$ \cK = R
  \oplus \La; $$ moreover, $R$ and $\La$ are both maximal isotropic
  subspaces of $\cK$.
\end{prop}

{\em Proof.} The fact that $\cK = R \oplus \La$ follows from
\cite{commI}, Prop.\ 1.1. That $\langle r_a, r_b \rangle =0$ follows from
\cite{examples}, 4.1.1. Let us show that $\langle r_a, \mm_1 \rangle$
vanishes: for $n>0$, $\on{val}_{P_1}(r_a z_1^n \omega) > (-n_1-1) +
n_1 = -1$ so $\res_{P_1}( r_a z_1^n \omega)$ is zero; and for $i>1$,
$\langle r_a, \cO_i \rangle$ vanishes because for $n\geq 0$,
$\on{val}_{P_i}(r_a z_i^n \omega) \geq -n_i + n_i = 0$ so $\res_{P_i}(
r_a z_i^n \omega)$ is zero.  \hfill \qed \medskip

\begin{remark} 
  In the case where $\omega$ has a unique zero of order $2(g-1)$ at
  some point $P_0$, $R$ is spanned by $f_0,f_{-a_1},f_{-a_2},...
  f_{-a_{g-1}}$, $f_{-2g},f_{-2g-1}, f_{-2g-2}, ...$ with
  $\on{val}_{P_0}(f_i) = i$: if $\omega_1,...,\omega_g$ be a basis of
  the space of holomorphic one-forms $H^0(X,\Omega_X)$, with
  $\on{val}_{P_0}(\omega_i) = b_i$, so that $0\leq b_1 < b_2< ... <
  b_g = 2(g-1)$, then $f_0 = 1$, $f_{-a_1} = \omega_{g-1}/\omega_g$,
  $f_{-a_2} = \omega_{g-1}/\omega_g$, etc.  On the other hand, the
  $r_a$ may be chosen to have poles of order $b_1, \ldots,b_g$ at
  $P_0$, with $\{a_1,\ldots,a_{g-1}\} \cup \{b_1,\ldots,b_g\} =
  \{1,\ldots,2g-1\}$.

  If $X$ is a hyperelliptic curve $y^2 = P_{2g+1}(x)$, $P_{2g+1}$ a
  polynomial of degree $2g+1$, and $\omega = dx/y$; more generally, if
  $X$ is a plane curve of equation $P(z) = Q(y)$, and $\omega =
  dx/Q'(y) = - dy/P'(x)$, with $P$ and $Q$ generic polynomials of
  coprime degrees $p$ and $q$, (in that case, $g = {{(p-1)(q-1)}\over
    2}$), $\omega$ has a zero of order $2(g-1)$ at the point at
  infinity.
\end{remark}

\subsection{(Twisted) Green functions}

We will denote by $z$ the $n$-uple $(z_i)$ of $\cK$. We will denote by
$\CC[[z,z^{-1}]]$ the set of series $\sum_{i=1}^p \sum_{n\in\ZZ}
a_{in} z_i^n$, and by $\CC[[z,w]]$ the space $\prod_{1 \leq i,j \leq
  p} \CC[[z_i,w_j]]$.
 
We define $\delta(z,w)$ as the sum $\sum_i \eps^i(z)\eps_i(w)$, where
$(\eps^i)$ and $(\eps_i)$ are dual bases of $\cK$ for $\langle ,
\rangle_\cK$.

The space of functions in two variables $z$ and $w$ will be identified
with the tensor square of the space of functions in one variable, via
the identification $a(z)b(w) \mapsto a\otimes b$.

\subsubsection{Green function}

Let $(e^i),(e_i)$ be dual bases of $R$ and $\La$. We will assume that
$(e_i)$ is the union of $(r_a)$ and a basis of $\mm = \mm_1 \oplus
\cO_2 \oplus \cdots \oplus \cO_p$. We set
\begin{equation} \label{G}
G  = \sum_i e^i \otimes e_i.  
\end{equation} 
We have then $\delta(z,w) = G(z,w) + G(w,z)$.  $G(z,w)$ is the
collection of expansions, for $w$ near each $P_i$, of a rational
function defined on $X^2$, antisymmetric in $z$ and $w$, regular
except for poles when $z$ of $w$ meets some $P_i$ and a simple pole at
the diagonal.

\subsubsection{Twisted Green functions}

To $\la_0 = (\la^{(0)}_a)_{1\leq a \leq g}$ a vector of $\CC^g$ is
associated the line bundle $\cL_{2\la_0}$ over $X$. The space $H^0(X-
\{P_0\}, \cL_{2\la_0})$ may be identified with the space of functions
on $\wt X$, regular outside $\pi^{-1}(P_0)$, with transformation
properties 
\begin{equation} \label{R:la}
  f(\gamma_{B_a}z) = f(z) \quad \on{and} \quad f(\gamma_{B_a}z) =
  e^{-2\la_a^{(0)}} f(z).
\end{equation} 
This space of functions will be denoted $R_{-2\la_0}$.  For $\la_0$
generic, a complement in $\cK$ of this space is $z_1^{1-g}\cO_1 \oplus
\cO_2 \oplus \cdots \oplus \cO_p$.

Let $\la = (\la_a)_{1\leq a \leq g}$ be $g$ formal parameters at the
vicinity of $\la_0$, and define $R_{-2\la}$ as the $\CC[[\la_a -
\la_a^{(0)}]]$-submodule of $\cK[[\la_a - \la_a^{(0)}]]$ generated by
the $e^{2\sum_a (\la_a - \la_a^{(0)}) r_a} \phi$, $\phi \in
R_{-2\la_0}$. Define also $\La' = (z_1^{1-g}\cO_1 \oplus \cO_2 \oplus
\cdots \oplus \cO_p)[[\la_a - \la_a^{(0)}]]$.

Then we have a direct sum decomposition 
$$ \cK[[\la_a - \la_a^{(0)}]] = R_{-2\la} \oplus \La'.
$$

For $\phi$ in $\cK[[\la_a - \la_a^{(0)}]]$, we denote by $\phi_{\La'}$
and $\phi_{R_{2\la}}$ the projections of $\phi$ on $\La'$ parallel to
$R_{2\la}$, resp.\ on $R_{2\la}$ parallel to $\La'$.  For $\phi(z)$ a
series $\sum_i \phi_i \eps_i(z)$, we define $\phi(z)_{z\to \La'}$ as
$\sum_i \phi_i (\eps_i)_{\La'}(z)$ and $\phi(z)_{z\to R_{2\la}}$ as
$\sum_i \phi_i (\eps_i)_{R_{2\la}}(z)$. 

We have then $f(z) = f(z)_{z\to \La'} + f(z)_{z\to R_{2\la}}$.

Let $(e^i_{2\la})$, $(e'_i)$ be dual bases of $R_{2\la}$ and $\La'$,
and let us set 
\begin{equation} \label{G:la}
G_{2\la}(z,w) = \sum_i e^i_{2\la}(z) e'_i(w). 
\end{equation}
We have $\delta(z,w) = G_{2\la}(z,w) + G_{-2\la}(w,z)$.
$G_{-2\la}(z,w) - G(z,w)$ belongs to $\CC[[z,w]][z^{-1},w^{-1}]$. 
The functions 
$$
g^+_{\la}(z) = (G_{-2\la} - G)(q^\pa z,z), \quad 
g^-_{\la}(z) = (G_{-2\la} - G)(q^{-\pa}z,z)
$$
then belong to $\cK[[\hbar]]$.

\begin{remark} 
  {\it Relation with the Green functions of \cite{commI}.} In
  \cite{commI}, we introduced Green function $G(z,w)$ and a twisted
  Green function $G_{2\la}(z,w)$, that we denote here by
  $G^{(\mathrm{I})}(z,w)$ and $G_{2\la}^{(\mathrm{I})}(z,w)$. Let
  $\omega_a$ be the basis of one-forms on $X$, associated to
  $(A_a)_{1\leq a \leq g}$.  The relation of these Green functions
  with $G(z,w)$ and $G_{2\la}(z,w)$ defined by (\ref{G}) (under the
  assumptions of Prop.\ \ref{Lambda}) and (\ref{G:la}) is
  $$ G(z,w) = \left( G^{(\mathrm{I})}(z,w) + \sum_a \omega_a(z)
    r_a(w)\right) / \omega(z)
  $$ and
  $$ G_{2\la}(z,w) = G_{2\la}^{(\mathrm{I})}(z,w)/\omega(z).$$
  
  Set $$\bar g_\la (z) = \lim_{z\to w}(G_\la(z,w) - G(z,w)).
  $$ One can show that
  $$ \bar g_{\la}(z) = \sum_a \pa_{\epsilon_a} \ln\Theta(-\la +
  (g-1)P_0- \Delta) \omega_a(z)/\omega(z).
  $$ where $(\eps_a)_{1\leq a \leq g}$ is the canonical basis of
  $\CC^g$, $\Theta$ is the Riemann theta-function on the Jacobian of
  $X$ and points of $X$ are identified with their images by the
  Abel-Jacobi map (\cite{commI}, 4.4).
\end{remark}

\section{The algebra $U_{\hbar,\omega}\G$} \label{pres:G}

\subsection*{Notation} 

For $E$ a vector space and $E',E''$ two subspaces, such that $E$ is
the direct sum $E' \oplus E''$, and for $\phi$ in $E$, we denote by
$\phi_{E'||E''}$ the projection of $\phi$ on $E'$ parallel to $E''$.
In the case of the decompositions $\cK = R \oplus \La$, $\cK[[\la_a -
\la_a^{(0)}]] = R_{2\la} \oplus \La'$ and $\cK = R_{(a)} \oplus \mm$
below, we will simply denote $\phi_{E'||E''}$ and $\phi_{E''||E'}$ by
$\phi_{E'}$ and $\phi_{E''}$ respectively.

For $(\eps_i)$ a basis of $\cK$ and $f(z)$ a series $\sum_i \eps_i(z)
\otimes v_i$ in some completion of $\cK \otimes V$, $V$ some vector
space, we define $f(z)_{z\to \La}$ as $\sum_i (\eps_i)_\La(z) \otimes
v_i$, and $f(z)_{z\to R}$ as $\sum_i (\eps_i)_R(z) \otimes v_i$. One
define in the same way $f(z)_{z\to R_{2\la}}$ and $f(z)_{z\to \La'}$.
If $f(z,w)$ is a series $\sum_{i,j} \eps_i(z)\eps_j(w)v_{ij}$,
$f(z,w)_{z\to\La}$ is $\sum_{i,j} (\eps_i)_{\La}(z)\eps_j(w)v_{ij}$,
etc.

If $f(z,w)$ belongs to $R_z((w))$ (the space of series $\sum_{i\geq
  n_0}r_i(z) w^i$, with $r_i$ in $R$) and there exists $g(z,w)$ in
$R_w((z))$, such that $f+g = (\pi\otimes id)\delta(z,w)$, where $\pi$
is some differential operator, then $f$ is the expansion of a rational
function on $X - \{P_i\})^2$ with only poles at the diagonal, and $g$
may be viewed as the analytic prolongation of $ - f$ in the region
$z<<w$. We write $g(z,w) = - f(z,w)_{z<<w}$. 

We write $\pa_z$ for $\pa \otimes id$, $\pa_w$ for $id \otimes \pa$.
we set $\phi^{(21)}(z,w) = \phi(w,z)$.

\subsection{Results on kernels} (see \cite{Ann3}) 
\label{kernels}

We have
$$
\pa_z G(z,w) = -G(z,w)^2 - \gamma, 
$$
for some $\gamma \in R\otimes R$. 
 
Let $\phi,\psi$ belong to $\hbar\CC[\gamma_0,\gamma_1,\ldots][[\hbar]]$
such that
$$
\pa_{\hbar}\psi = D\psi - 1 - \gamma_0 \psi^2, \quad
\pa_{\hbar}\phi = D\phi - \gamma_0 \psi. 
$$
Here $D = \sum_{i\ge 0}\gamma_{i+1} \pa_{\gamma_i}$.  
we have 
$$ \psi(\hbar,\pa_z^i \gamma) = -\hbar +o(\hbar), \quad
\phi(\hbar,\pa_z^i\gamma) = {1\over 2}\hbar^2 \gamma_0 + o(\hbar^2).
$$

Set $G^{(21)}(z,w) = G(w,z)$. 
%We have
%$$ \sum_i {{1-q^{-\pa}}\over{\pa}}e_i(z) \otimes e^i(w) = -
%\phi(\hbar, (-\pa_z)^i \gamma) + \ln(1 - G^{(21)}\psi(\hbar,
%(-\pa_z)^i\gamma)),
%$$
From identity (3.11) of \cite{Ann3} (with $\pa$ transformed to $-\pa$)
and by (3.8) of \cite{Ann3}, we have
\begin{equation} \label{giov} 
\sum_i {{1-q^{-\pa}}\over{\pa}}e_i(z) \otimes e^i(w) = -
\phi( - \hbar, \pa_z^i \gamma) + \ln(1 + G^{(21)}\psi(-\hbar,
\pa_z^i\gamma)). 
\end{equation}

Set $T = {{\mathrm{sinh} \hbar\pa}\over{\hbar\pa}}$.  Let $\tau$ in
$(R\otimes R)[[\hbar]]$ satisfy
\begin{equation} \label{cond:tau} \label{id:U}
\tau + \tau^{(21)} = - \sum_i e^i \otimes (Te_i)_R. 
\end{equation} 
Let $U$ be the linear map from $\La$ to $R[[\hbar]]$ such that $\tau =
\sum_i Ue_i \otimes e^i$. We have
$$
\sum_{i} (T+U)e_i \otimes e^i  + \sum_i e^i \otimes (T+U)e_i = 
(T\otimes id)\delta(z,w), 
$$ which means that after analytic prolongation the sum $\sum_{i}
(T+U)e_i \otimes e^i$ is antisymmetric in $z$ and $w$.

Set $T_+ = {{1-q^{-\pa}}\over{2\hbar\pa}}$ and define $U_+ : \La\to
R[[\hbar]]$ by the formula $U_+ = (1+q^\pa)^{-1} \circ U$; we have 
$$ (T_+ + U_+)(\la) = {1\over{1+q^\pa}} ((T+U)(\la)). $$ Define $q_+$
by
$$
q_+(z,w) = q^{2\sum_i (T_+ + U_+)e_i(z) \otimes e^i(w)}, 
$$ it then follows from (\ref{giov}) that
\begin{equation} \label{id:q+} 
q_+(z,w) = q^{2\sum_i (U_+e_i)(z) \otimes e^i(w)} e^{-\phi(- \hbar,
  \pa_z^i \gamma)}(1 + G^{(21)}\psi(-\hbar, \pa_z^i\gamma)).
\end{equation}

\begin{remark} 
  Formulas of this section correct a sign mistake in \cite{Ann3}: in
  sect.\ 3 of that paper, $\pa$ should be changed to $-\pa$.
\end{remark}

\subsection{} 

The algebra $U_{\hbar,\omega}\G$ has generators
$h^+[r],h^-[\la],e[\eps],f[\eps]$ and $K$, with $r$ in $R$, $\la$ in
$\La$ and $\eps$ in $\cK$; generating series 
$$
x(z) = \sum_i x[\eps^i]\eps_i(z),\ h^+(z) = \sum_i h^+[e^i]e_i(z), \ 
h^-(z) = \sum_i h^-[e_i]e^i(z),  
$$ $x = e,f$, and relations
$$
x[\al\eps + \eps'] = \al x[\eps] + x[\eps'], 
$$
for $\al$ scalar, $x = h^+$, $\eps,\eps'$ in $R$; 
$x = h^-$, $\eps,\eps'$ in $\La$; or $x = e,f$, $\eps,\eps'$ in $\cK$; 
\begin{equation} \label{comm}
  [h^{+}[r], h^{+}[r']]=0,
\end{equation}
\begin{equation} \label{h-h}
  [K,\on{anything}]=0, \quad [h^{+}[r], h^{-}[\la]]= {2 \over \hbar}
  \langle (1 - q^{-K \partial })r, \la \rangle,
 \end{equation}
\begin{equation} \label{h-h-}
  [h^{-}[\la], h^{-}[\la']]={2\over\hbar} \left( \langle
    T((q^{K\pa}\la)_R) , q^{K\pa}\la' \rangle + \langle
    U\la,\la'\rangle - \langle U((q^{K\pa}\la)_\La),
    q^{K\pa}\la'\rangle \right)
\end{equation}
\begin{equation}   \label{h-e}
  [h^{+}[r], e(w)]=2 r(w)e(w), \quad [h^{-}[\la],
  e(w)]=2[(T+U)(q^{K\pa}\la)_{\La}](w)e(w),
\end{equation}
\begin{equation}   \label{h-f}
  [h^{+}[r], f(w)]=-2 r(w)f(w), \quad [h^{-}[\la],
  f(w)]=-2[(T+U)\la](w)f(w),
\end{equation}
\begin{equation}  \label{e-e}
(\al(z) - \al(q^{-\pa}w)) e(z)e(w) = (\al(z) - \al(q^{-\pa}w))
q^{2\sum_i (T+U)e_i(z)\otimes e^i(w)}e(w)e(z)
\end{equation}
and
\begin{equation}  \label{f-f}
(\al(z) - \al(q^\pa w)) f(z)f(w) = (\al(z) - \al(q^\pa
w))q^{-2\sum_i (T+U)e_i(z)\otimes e^i(w)}f(w)f(z)
\end{equation}
for any $\al$ in $\cK$, 
\begin{equation}   \label{e-f}
  [e(z),f(w)]={1\over \hbar} [ \delta(z,w)q^{((T+U)h^{+})(z)}
  -(q^{-K\pa_{w}}\delta(z,w))q^{-h^{-}(w)} ] ,  
\end{equation}
$\phi(z,w) = \phi(\hbar,\pa_z^i\gamma)$, $r,r'$ in $R$, $\la,\la'$ in
$\La$ (see \cite{Ann3}).  $U_{\hbar,\omega}\G$ is completed with
respect to the topology defined by the left ideals generated by the
$x[\eps]$, $\eps$ in $\oplus_i z_i^N\cO_i$ The critical case
correponds to $K=-2$.

\begin{remark}
Relations (\ref{e-e}) and (\ref{f-f}) can be written 

\end{remark}

\subsection{Cartan currents}

In case we have a relation
$$ a(z) b(w) a(z)^{-1} = \mu(z,w) b(w),
$$ with $a(z),b(w)$ currents of $U_{\hbar,\omega}\G$ and $\mu(z,w)$ in
$\CC((z))((w))[[\hbar]]$ or $\CC((w))((z))[[\hbar]]$, we will define 
$( a(z), b(w))$ as $\mu(z,w)$.  

Set $K^+(z) = q^{(T+U)h^+(z)} , K^-(z) = q^{-h^-(z)}$. Let us also set
$$ q(z,w) = q^{2\sum_i (T+U)e_i(z)\otimes e^i(w)} ;
$$ we have $q(z,w) = (q(w,z)^{-1})_{w<<z}$.  Then the
relations involving Cartan generators can be expressed as
\begin{equation} \label{K-K}
  (K^{+}(z),K^{+}(w))=1, \quad (K^{+}(z),K^{-}(w))={{q(z,q^{-K\pa}(w))}
    \over {q(z,w)}},
\end{equation}
\begin{equation} \label{K-K-}
  (K^{-}(z),K^{-}(w))={{q(q^{-K\pa}(z),q^{-K\pa}(w))}\over{q(z,w)}},
\end{equation}
\begin{equation}\label{K-e}
  (K^{+}(z),e(w))= q(z,w), \quad (K^{-}(z),e(w))=
  q(w,q^{-K\pa}(z))^{-1} ,
\end{equation}
\begin{equation}\label{K-f}
  (K^{+}(z),f(w))= q(z,w)^{-1}, \quad (K^{-}(z),f(w))= q(w,z). 
\end{equation}

Set
$$ k^+(z) = q^{(T_+ + U_+)h^+(z)}, \quad k^-(z) = \la(z)
q^{{1\over{1+q^{-\pa}}} h^-(z)}, 
$$ 
with $\la(z)$ the function such that 
$$
\la(z) \la(q^{-\pa}z) q^{[{1\over{1+q^{-\pa}}} h^-(z) ,
{{q^{-\pa}}\over{1+q^{-\pa}}} h^-(z)]} = 1, 
$$
that is
$$ \la(z) = \exp\left[- {1\over{1+q^{-\pa_z}}}\left(
    ({1\over{1+q^{-\pa_z}}}\otimes
    {{q^{-\pa_{z'}}}\over{1+q^{-\pa_{z'}}}})[h^-(z),h^-(z')]\right)_{z'
    = z} \right]. 
$$

We have
$$ K^+(z) = k^+(z)k^+(q^\pa z), \quad K^-(z) = k^-(z)^{-1}k^-(q^{-\pa}
z)^{-1}.
$$

Set
$$ q_+(z,w) = q^{2\sum_i (T_+ + U_+)e_i(z) \otimes e^i(w)}, q_-(z,w) =
q^{-2\sum_i {1\over{1+q^{\pa}}}e^i(z) \otimes (T + U)e_i(w)},
$$ then we have $q_+(z,w) = q_-(z,w)$ (up to analytic continuation).
We have
$$ q_+(z,w)q_+(q^{\pa}z,w) = q(z,w),
$$ and
$$ (k^+(z),e(w)) = q_+(z,w), \quad (k^-(z)^{-1},e(w)) =
q_-(q^{(-K+1)\pa}z,w),
$$
$$ (k^+(z),f(w)) = q_+(z,w)^{-1}, \quad (k^-(z)^{-1},f(w)) =
q_-(q^{\pa}z,w)^{-1}.
$$

Also when $K = -2$, we have 
\begin{equation} \label{17b} 
(k^+(z), k^+(w)) = 1, \quad 
(k^+(z),k^-(w)) = {{q_+(z,q^\pa w)}\over{q_+(z,q^{2\pa}w)}},
\end{equation}
and
\begin{equation} \label{17c} 
(k^-(z),k^-(w)) = {{q_+(q^{3\pa}z,q^{2\pa}w)} \over
  {q_+(q^{\pa}z,q^{\pa}w)}}
{{q_-(q^{2\pa}w,q^{2\pa}z)}\over{q_-(q^{2\pa}w,q^\pa z)}} . 
\end{equation}

% \quad
%(k^-(z), K^-(w)) = {{q_+(q^\pa z,w)}\over{q_+(q^{3\pa}z, q^{2\pa}w)}}.
%$$

\section{Subalgebra $U_\hbar\B_{in}$} \label{sub+}

The quantity
$$
(q^{\pa_z + \pa_w} - 1) \sum_i ((T+U)e_i)(z) e^i(w)
$$ belongs to $(R\otimes R)[[\hbar]]$, since $Ue_i$ belongs to $R$,
$T$ commutes with $\pa_z + \pa_w$ and $(\pa_z + \pa_w)G(w,z) = -
\sum_i e^i(z) (\pa e_i)_R(w)$ belongs to $R\otimes R$. Moreover 
$$ F(z,w) = 2\hbar {{q^{\pa_z + \pa_w} -
    1}\over{(1+q^{-\pa_z})(1+q^{-\pa_w})}} \sum_i (T+U)e_i (z) e^i(w)
$$
is symmetric in $z$ and $w$. 
Let $\al(z,w)$ be an element of $\hbar(R\otimes R)[[\hbar]]$ such that
\begin{equation} \label{def:al} 
  {{\exp(2\al(q^\pa w, z)) }\over{\exp(2\al(q^\pa z, w))}} =
  \exp[2\hbar {{q^{\pa_z + \pa_w} -
      1}\over{(1+q^{-\pa_z})(1+q^{-\pa_w})}} \sum_i (T+U)e_i (z)
  e^i(w)] ; 
\end{equation}
we may choose 
$$ \al(z,w) = {1\over 2} F(w,q^{-\pa}z).
$$

Let us set
$$ \al(z,w) = \sum_{i,j} a_{ij} e^i(z) e^j(w), 
$$
and
$$ k_R(z) = \exp(\sum_{i,j} a_{ij}h[e^i] e^j(z)) .
$$

Define $R_{(a)}$ as $\oplus_a \CC r_a \oplus R$.  Recall that we
defined $\mm$ as $\mm_1 \oplus \cO_2 \oplus \cdots \oplus \cO_p$, so
that $\cK = R_{(a)} \oplus \mm$. Let $\cA$ be the
$\CC[[\hbar]]$-module automorphism of $R_{(a)}[[\hbar]]$ defined by
\begin{equation} \label{cA}
  \cA(r) = r\ \on{ for} \ r \ \on{in} \ R, \ \on{and} \ \cA(r_a) =
  (T+U)r_a\ \on{ for} a = 1, \ldots,g.
\end{equation}

Define $\beta(z,w)$ in $(R \otimes R_{(a)})[[\hbar]]$ by 
\begin{align*} 
  & -2(1 \otimes q^\pa \cA) \beta(z,w) \\ & = 2(\al(q^{2\pa}z,w) -
  \al(q^\pa z,w)) - 2 \hbar \sum_i {1\over{1+q^{-\pa}}} e^i(z) \otimes
  ((T+U)e_i)_{R_{(a)}}(w) \\ & + 2 \hbar (q^{3\pa} \otimes q^\pa - q^\pa
  \otimes q^{-\pa})({1\over{1 + q^\pa}} \otimes
  {1\over{1+q^{-\pa}}})\sum_i ((T+U)e_i)(z) \otimes e^i(w).
\end{align*}
Set 
$$ \beta(z,w) = \sum_{a,i} b_{ai} e^i(z) r_a(w) + \sum_{i,j} c_{ij}
e^j(z) e^i(w)
$$ and
\begin{equation} \label{def:ka} 
k_a(z) = \exp(\sum_{a,i} b_{ai}h[r_a]e^i(z) + \sum_{i,j} c_{ij}
h[e^i] e^j(z)) . 
\end{equation}
%we have then 
%\begin{align} \label{def:la} 
%  & (k_a(z), f(w)k^-(q^{-\pa}w)) = \exp[2(\al(q^{2\pa}z,w)-
%  \al(q^{\pa}z,w))] \cdot \\ & \nonumber \cdot \exp[-2\hbar \sum_i
%  {1\over{1+q^{-\pa}}} e^i(z) \otimes ((T+U)e_i)_{R_{(a)}}(w)] (k^-(z),
%  k^-(q^{-\pa}w)).
%\end{align}

Set finally 
$$
k_{\mm}(z) = k_a(z)^{-1} k^-(z). 
$$ The currents $k_R(z),k_a(z)$ and $k_{\mm}(z)$ all belong to
$U_\hbar\HH \otimes R_z[[\hbar]]$.

\begin{prop} \label{walsingham} 
  i) Set $\wt f(z) = f(z) k_R(z) k^-(q^{-\pa}z)$. we have $\wt f(z)
  \wt f(w) = \wt f(w) \wt f(z)$.

  ii) We have
  $$ (k_{\mm}(z), \wt f(w)) \in \exp(\hbar (R\otimes\mm) [[\hbar]]).
  $$ 

  iii) We have $k_\mm(z)k_\mm(w) = k_\mm(w)k_\mm(z) $.
\end{prop}

{\em Proof.} Let us show that $\wt f(z)$ commutes with itself.  We
have
$$[h[r] , f(z) k^-(q^{-\pa}z)] = -2 (q^\pa r)(z)
f(z)k^-(q^{-\pa}z),  
$$ 
therefore 
$$ { {(k_R(w),f(z) k^-(q^{-\pa}z))} \over {(k_R(z),f(w)
    k^-(q^{-\pa}w))} } = {{\exp(2\al(q^\pa z,w))}\over{\exp(2\al(q^\pa
    w,z))}}.
$$ Set
\begin{equation} \label{j}
j(z,w) = q_+(q^\pa z,w)q_-(w,z), 
\end{equation} 
we have $j(z,w)\in 1 + \hbar (R\otimes R)[[\hbar]]$, $j(z,w) j(w,z) =
1$.  From \cite{examples} follows that for some $i_-(z,w)$ in
$\CC[[z,w]][z^{-1},w^{-1}][[\hbar]]^\times$, we have
$$
q_-(z,w) = i_-(z,w) {{w - q^{-\pa}z}\over{w - z}}, 
$$ 
so that 
$$ j(z,w) = i_-(q^\pa z, w) i_-(w,z){{z - q^{-\pa}w}\over{q^\pa z -
    w}}. 
$$
We have then 
\begin{align*}
  & (w-z) f(z) k^-(q^{-\pa}z) f(w) k^-(q^{-\pa}w) \\ & = i_-(z,w) (w -
  q^{-\pa}z) f(z) f(w) k^-(q^{-\pa}z)k^-(q^{-\pa}w) \\ & = (
  k^-(q^{-\pa}z), k^-(q^{-\pa}w)) i_-(q^\pa z,w)^{-1} (w - q^\pa z)
  f(w) f(z) k^-(q^{-\pa}w)k^-(q^{-\pa}z) \\ & = ( k^-(q^{-\pa}z),
  k^-(q^{-\pa}w)) i_-(q^\pa z,w)^{-1} (w - z) {{w - q^\pa
      z}\over{q^{-\pa}w - z}}f(w) k^-(q^{-\pa}w) f(z) k^-(q^{-\pa}z), 
\end{align*}
therefore
\begin{align*} 
  (w-z) & \left[ \left( f(z)k^-(q^{-\pa}z) \right) \left(
      f(w)k^-(q^{-\pa}w) \right) \right. \\ & \left. - {{(
        k^-(q^{-\pa}z) , k^-(q^{-\pa}w))}\over{j(z,w)}} \left(
      f(w)k^-(q^{-\pa}w) \right) \left( f(z)k^-(q^{-\pa}z)
    \right)\right] = 0;
\end{align*}
let $B(z,w)$ be the term in brackets. It is equal to
$A(z)\delta(z,w)$, for some generating series $A(z)$. Since $B(z,w)$
also satisfies $B(w,z) = - {{( k^-(q^{-\pa}w), k^-(q^{-\pa}z)
    )}\over{j(w,z)}} B(z,w)$, and ${{( k^-(q^{-\pa}w), k^-(q^{-\pa}z)
    )}\over{j(w,z)}} = 1 + o(\hbar)$, we obtain that $B(z,w)$ vanishes
so
\begin{align} \label{starr}
  & \left( f(z)k^-(q^{-\pa}z) \right) \left( f(w)k^-(q^{-\pa}w) \right)
  \\ & \nonumber = {{( k^-(q^{-\pa}z) , k^-(q^{-\pa}w))}\over{j(z,w)}} \left(
    f(w)k^-(q^{-\pa}w) \right) \left( f(z)k^-(q^{-\pa}z) \right) .
\end{align}
We have
$$ (k^-(z), k^-(w)) = \exp[2\hbar (q^{2(\pa_z + \pa_w)} - 1)
{{1}\over{1+q^{-\pa_z}}} {{1}\over{1+q^{-\pa_w}}} \sum_i (T+U)e_i(z)
e^i(w)], 
$$
\begin{equation} \label{id:j} 
  j(z,w) = \exp[2\hbar (1- q^{- \pa_z - \pa_w})
  {{1}\over{1+q^{-\pa_z}}} {{1}\over{1+q^{-\pa_w}}} \sum_i (T+U)e_i(z)
  e^i(w)],
\end{equation}
so
$$ {{( k^-(q^{-\pa}z), k^-(q^{-\pa}w))}\over{j(z,w)}} = \exp[2\hbar
{{q^{\pa_z + \pa_w} - 1}\over{ (1+q^{-\pa_z})(1+q^{-\pa_w})}} \sum_i
(T+U)e_i(z) e^i(w)]; 
$$ since $k_R(z)$ commutes with $k_R(w)$, and by (\ref{starr}), {\em
  i)} follows.

Let us prove {\em ii)}. We have 
$$ (k^-(z) , \wt f(w) ) = \exp[2(\al(q^{2\pa}z,w)- \al(q^{\pa}z,w))]
q_-(q^\pa z , w) (k^-(z) , k^-(q^{-\pa}w)) ; 
$$ moreover, 
$$ q_-(q^\pa z, w) = \exp[-2\hbar \sum_i {1\over{1 + q^{-\pa}}} e^i(z)
\otimes (T+U)e_i(w) ],
$$ From (\ref{def:ka}) follows that
\begin{align*} & (k_a(z), \wt f(w)) = \exp[2(\al(q^{2\pa}z,w)-
  \al(q^{\pa}z,w))] \cdot \\ & \cdot \exp[-2\hbar \sum_i
  {1\over{1+q^{-\pa}}} e^i(z) \otimes ((T+U)e_i)_{R_{(a)}}(w)] (k^-(z),
  k^-(q^{-\pa}w)).
\end{align*}
Therefore,
\begin{equation} \label{ccc}
(k_\mm(z), \wt f(w)) = \exp[-2\hbar \sum_i {1\over{1+q^{-\pa}}} e^i(z) 
\otimes ((T+U)e_i)_\mm(w)], 
\end{equation} which implies {\em ii)}.

Set for $\phi = \sum_a \la_a r_a + r$, with $\la_a$ in $\CC$ and $r$
in $R$, $h[\phi] = \sum_a \la_a h^-[r_a] + h^+[r]$.  Then we have for
$\phi$ in $R_{(a)}$,
$$ [h[\phi], f(z) ] = -2 (\cA\phi)(z) f(z)
$$ 
and 
$$ [h[\phi], k^-(z)] = 2 [q^\pa (1 - q^\pa) \cA\phi](z) k^-(z),
$$
where $\cA$ is defined by (\ref{cA}), 
so that 
$$ [h[\phi], \wt f(z)] = -2 [q^\pa \cA\phi](z) \wt f(z).
$$ Therefore we get
\begin{align*} 
  (k_a(z) , k^-(w)) & ={{(k_a(z), \wt f(q^\pa w))}\over{(k_a(z), \wt
      f(w))}} \\ & = \exp[2(q^{2\pa_z + \pa_w}- q^{\pa_z + \pa_w} -
  q^{2\pa_z} + q^{\pa_z}) \al(z,w)] \cdot \\ & \cdot \exp[-2\hbar
  \sum_i {1\over{1+ q^{-\pa}}} e^i(z) \otimes (q^\pa - 1)
  ((T+U)e_i)_{R_{(a)}}(w)] \cdot \\ & \cdot
  {{(k^-(z),k^-(w))}\over{(k^-(z),k^-(q^{-\pa}w))}}. 
\end{align*}
On the other hand, {\em iii)} is translated as
$$
{{(k_a(z)^{-1}, k^-(w))} \over {(k_a(w)^{-1}, k^-(z))}} 
(k^-(z), k^-(w)) = 1,  
$$ that is
 \begin{align*}
   & \exp[(q^{\pa_z} - 1)(q^{\pa_w} - 1)\left( 2\al(q^\pa z,w)-
     2\al(q^\pa w,z)\right) ] \\ & \exp[-2\hbar \sum_i {1\over{1+
       q^{-\pa}}} e^i(z) \otimes (q^\pa - 1) ((T+U)e_i)_{R_{(a)}}(w)] :
   (z\leftrightarrow w) \\ &
   {{(k^-(z),k^-(w))^2}\over{(k^-(z),k^-(q^{-\pa }w)) (k^-(q^{-\pa }z)
       , k^-(w))}} \\ & = (k^-(z),k^-(w)),
\end{align*}
in other terms 
\begin{align*}
  & \exp[(q^{\pa_z} - 1)(q^{\pa_w} - 1) \log \left( {{(k^-(q^{-\pa}z),
        k^-(q^{-\pa}w))}\over{j(z,w)}} \right)^{-1} ] \\ &
  \exp[-2\hbar \sum_i {1\over{1+ q^{-\pa}}} e^i(z) \otimes (q^\pa - 1)
  ((T+U)e_i)_{R_{(a)}}(w)] : (z\leftrightarrow w) \\ &
  {{(k^-(z),k^-(w))}\over{(k^-(z),k^-(q^{-\pa }w)) (k^-(q^{-\pa }z) ,
      k^-(w))}} \\ & = 1,
\end{align*}
or
\begin{align} \label{3'}
  & \exp[-2\hbar\sum_i {1\over{1+ q^{-\pa}}} e^i \otimes (q^\pa
    -1)((T+U)e_i)_{R_{(a)}}] : (z \leftrightarrow w)
\\ & \nonumber \exp[(q^{\pa_z}-1)(q^{\pa_w}-1) \log j(z,w)]
\\ & \nonumber = ( k^-(q^{-\pa}z) , k^-(q^{-\pa}w)). 
\end{align}

The terms containing $U$ in the logarithm of (\ref{3'}) are
\begin{align*}
  & -2 \hbar ({1\over{1 + q^{-\pa}}} \otimes (q^\pa - 1)) \sum_i e^i
  \otimes Ue_i \\ & + 2 \hbar ((q^\pa - 1) \otimes
  {1\over{1+q^{-\pa}}}) \sum_i Ue_i \otimes e^i + 2 \hbar ({{q^\pa -
      1}\over{1 + q^{-\pa}}} \otimes (q^\pa - 1))\sum_i Ue_i \otimes
  e^i \\ & - 2 \hbar ((q^\pa - 1)\otimes{{q^\pa - 1} \over{q^\pa +
      1}}) \sum_i Ue_i \otimes e^i \\ & -2\hbar (q^{\pa} \otimes
  q^{\pa}- q^{-\pa} \otimes q^{-\pa}) ({1\over {1 + q^{-\pa}}}\otimes
  {1\over {1 + q^{-\pa}}}) (\sum_i Ue_i \otimes e^i)
\end{align*}
which is equal to 
$$ 2\hbar ({1\over{1+q^{-\pa}}} \otimes(q^\pa - 1)) \sum_i e^i \otimes
(Te_i)_R, 
$$ in view of (\ref{id:U}). 

Therefore {\em iii)} is written as
\begin{align*}
  & \exp[-2\hbar ({1\over{1 + q^{-\pa}}} \otimes (q^\pa - 1)) \sum_i
  e^i \otimes (Te_i)_{R_{(a)}}] : (z\leftrightarrow w) \\ & \exp[2\hbar
  ({1\over{1+q^{-\pa}}}\otimes(q^\pa - 1)) \sum_i e^i \otimes
  (Te_i)_R] \\ & \exp[2\hbar ({{q^\pa - 1 }\over{1+ q^{-\pa}}}
  \otimes(q^\pa - 1)) \sum_i Te_i \otimes e^i] \exp[ - 2\hbar
  ((q^\pa - 1) \otimes {{q^\pa - 1 }\over{1+q^\pa}}) \sum_i Te_i
  \otimes e^i] \\ & = \exp[2\hbar (q^\pa \otimes q^\pa - q^{-\pa}
  \otimes q^{-\pa}) ({1\over{1+q^{-\pa}}}\otimes
  {1\over{1+q^{-\pa}}})\sum_i Te_i \otimes e^i]
\end{align*}
or
\begin{align} \label{stg}
  & \nonumber \exp[2\hbar ({1\over{1+q^{-\pa}}}\otimes (1 - q^\pa))
  (\sum_i Te_i \otimes e^i + e^i \otimes (Te_i)_{R_{(a)}} - e^i \otimes
  (Te_i)_R) ] \\ & \exp[2\hbar ((1 - q^\pa)\otimes
  {1\over{1+q^{-\pa}}}) (\sum_i Te_i \otimes e^i - (Te_i)_{R_{(a)}}
  \otimes e^i ) ] = 1.
\end{align}
Since 
$$ \sum_i Te_i \otimes e^i + e^i \otimes (Te_i)_{R_{(a)}} - e^i \otimes
(Te_i)_R
$$
is
\begin{equation} \label{contrapunto} 
(T\otimes id)\delta(z,w) - \left( \sum_i Te_i \otimes e^i -
  (Te_i)_{R_{(a)}} \otimes e^i \right)^{(21)}, 
\end{equation}
(\ref{stg}) is equivalent to the statement that
\begin{align} \label{key}
  & ((1 - q^\pa)\otimes {1\over{1+q^{-\pa}}}) \left( \sum_i Te_i
    \otimes e^i - (Te_i)_{R_{(a)}} \otimes e^i \right) -
  (z\leftrightarrow w) \\ & \nonumber + ({{(q^\pa - 1)(1 -
      q^{-\pa})}\over{2\hbar \pa}} \otimes id) \delta(z,w) = 0
\end{align}
(whose interpretation is that after analytic prolongation, 
$$ ((1 - q^\pa)\otimes {1\over{1+q^{-\pa}}}) \left( \sum_i Te_i
  \otimes e^i - (Te_i)_{R_{(a)}} \otimes e^i \right)
$$ is symmetric in $z$ and $w$). 

To prove this, we first show that 

\begin{lemma} \label{pp}
  One can choose the dual bases $(e_i)_{i\geq 0}$, $(e^i)_{i\geq 0}$
  as $(r_a;e'_i)_{a = 1,...,g, i\geq 0}$ and $(\omega_a / \omega;
  e^{\prime i})_{a = 1,...,g, i\geq 0}$, with $(e'_i)_{i\geq 0}$ be a
  basis of $\mm$ and $(e^{\prime i})_{i\geq 0}$ the dual basis of the
  subspace $K_a$ of $R$ defined as $\{r\in R | \int_{A_i} r\omega =
  0\}$. We have
\begin{equation} \label{equality}
  \sum_{i\geq 0} (Te_i)_{\mm} \otimes e^i = \sum_{i\geq 0} e'_i
  \otimes Te^{\prime i} .
\end{equation}
\end{lemma}

{\em Proof of Lemma.} The first statement follows from the fact that
$K_a$ is the annihilator of $\oplus_a \CC r_a$ in $R$ for $\langle ,
\rangle_\cK$. Let us show the second statement. Both sides of the
equality belong to $\mm \otimes \cK$.  On the other hand, the
annihilator of $K_a$ for $\langle , \rangle_{\cK}$ is $R_{(a)}$ and has
therefore zero intersection with $\mm$. It follows that to show
(\ref{equality}), it is enough to show that the pairing of both sides
with $\rho \otimes id$ coincide, for $\rho$ in $K_a$. But $\langle
\on{left} \ \on{ side}, \rho \otimes id\rangle = \sum_i (T e_i)_\mm ,
\rho \rangle e^i = \sum_i \langle Te_i, \rho\rangle e^i$, because
$K_a$ and $R_{(a)}$ are orthogonal; this is equal to $\sum_i \langle e_i,
T\rho\rangle e^i$ because $T$ is self-adjoint and therefore to
$T\rho$.  On the other hand, $\langle \on{right}\ \on{side}, \rho
\otimes id \rangle = \sum_i \langle e'_i,\rho \rangle Te^{\prime i} =
T\rho$. This proves (\ref{equality}).
%Since $Tr_a$ belongs to $R_{(a)}$, belongs to $\mm \otimes K_a$.  
\hfill \qed \medskip
 
(\ref{key}) is equal to
$$ ((1 - q^\pa)\otimes {1\over{1+q^{-\pa}}}) \sum_i (Te_i)_{\mm}
\otimes e^i ;
$$
by Lemma \ref{pp}, this is 
$$ ((1 - q^\pa)\otimes {1\over{1+q^{-\pa}}}) \sum_i e'_i \otimes
Te^{\prime i}, 
$$
which is 
$$ (1 - q^\pa) \otimes {{q^\pa - 1}\over{2\hbar\pa}} \sum_i e'_i
\otimes e^{\prime i}
$$
or
$$ - {1\over {2\hbar}} ({{q^\pa - 1} \over \pa} \otimes {{q^\pa -
    1}\over{\pa}}) \sum_i \pa e'_i \otimes e^{\prime i}. 
$$

(\ref{key}) now follows from
$$ \sum_i \pa e'_i \otimes e^{\prime i} -
\sum_i e^{\prime i} \otimes \pa e'_i = (\pa\otimes id) \delta(z,w)
$$ (which means that after analytic continuation, $\sum_i \pa e'_i
\otimes e^{\prime i}$ is symmetric). This equality can be proved
either by expressing $\sum_i e'_i \otimes e^{\prime i}$ explicitly
using theta-functions (see \cite{commI}), or as follows: $\sum_i \pa
e'_i \otimes e^{\prime i} - e^{\prime i}\otimes \pa e'_i -
(\pa\otimes id) \delta(z,w)$ is equal to
\begin{equation} \label{balanus}
  \sum_i \pa e'_i \otimes e^{\prime i} + e'_i\otimes \pa e^{\prime i}
  + \sum_a \omega_a / \omega \otimes \pa r_a,
\end{equation} 
which belongs to $\cK \otimes R$. To show that (\ref{balanus}) is zero,
let us pair in with $id \otimes \la$, $\la$ in $\mm \oplus \oplus_a \CC
r_a$. For $o$ in $\mm$, $\langle (\ref{balanus}), id \otimes o \rangle$
is equal to
$$ (\pa o)_{R_{(a)} || \mm} - \sum_a \omega_a / \omega \langle o , \pa r_a
\rangle ; 
$$ 
 but $\pa o$ belongs to $\mm \oplus (\oplus_a \CC r_a) \oplus
(\oplus_a \CC \omega_a / \omega)$, therefore this vanishes.  On the
other hand, $\langle (\ref{balanus}), id \otimes r_a \rangle$ is equal
to zero, because $\langle \pa r_a , r_b \rangle =0$ for any $a,b$.  \hfill
\qed \medskip

{}From the proof of Prop.\ \ref{walsingham} follows that
$(k_\mm(z),\wt f(w))$ is of the form $\exp(-2\hbar \sum_i e^{\prime i}
\otimes e'_i + o(\hbar))$, so that we have
\begin{equation} \label{kO:wtf} 
(k_\mm(z), \wt f(w)) = \exp(-2\hbar \sum_{i\geq 0} \varphi^i
\otimes e'_i),
\end{equation}
with $(\varphi^i)_{i\geq 0}$ a free family of $R[[\hbar]]$. Set
$h_\mm(z) = {1\over\hbar} \ln k_\mm(z)$; (\ref{kO:wtf}) implies that
$$ h_\mm(z) = \sum_{i\geq 0} \wt h[e'_i] \varphi^i(z),
$$ with $\wt h[e'_i]$ linear combinations of the $h^+[r]$ and
$h^-[\la]$. Define $\wt h[o]$ for $o$ in $\mm$ by linear extension.

\begin{cor}
  Define in $U_{\hbar,\omega}\G$, $\wt f[\eps]$ as $\sum_i
  \res_{P_i}(\wt f(z) \eps(z) \omega(z))$. Let $\B_{in}$ be the Lie
  algebra $$\B_{in} = (\bar\HH \otimes \mm) \oplus (\bar\N_+ \otimes
  \cK).$$ Then there is an algebra injection $U\B_{in}[[\hbar]] \to
  U_{\hbar,\omega} \G$ defined by $h\otimes o \mapsto \wt h[o],
  f\otimes \eps \mapsto \wt f[\eps]$. We define $U_\hbar\B_{in}$ as
  the image of this injection.
\end{cor}

{\em Proof.} From the construction of $\wt h[o]$ follows that we have
$$
[\wt h[o],\wt f(z)] = -2 o(z) \wt f(z);  
$$ Prop.\ \ref{walsingham}, {\em i)} and {\em iii)} then imply the
statement.  \hfill \qed \medskip

Define $U_\hbar\B_-$ as the subalgebra of $U_{\hbar,\omega}\G$
generated by the $h^+[r],h^-[\la]$ and the $f[\eps]$, $r$ in $R$,
$\la$ in $\La$, $\eps$ in $\cK$; $U_\hbar\B_{in}$ is then a subalgebra
of $U_{\hbar}\B_-$. Define $U_\hbar \B_{\la_0}^{out}$ as the
subalgebra of $U_\hbar\B_-$ generated by the $h^+[r]$ and the $\wt
f[r_{-2\la_0}]$, $r$ in $R$, $r_{-2\la_0}$ in $R_{-2\la_0}$.
%$k_{\mm}(z)$ and $\wt f(z)_{z\to  R_{2\la}}$  

We have:

\begin{prop} \label{PBW:borel} 
  $U_\hbar\B_-$ is the direct sum of $\CC[h[r_a]] U_\hbar\B^{in}$ and
  of its right ideal generated by its right ideal $\sum_{r\in R}h^+[r]
  U_\hbar \B^{out}_{\la_0} + \sum_{r_{-2\la_0} \in R_{-2\la_0}} \wt
  f[r_{-2\la_0}] U_\hbar \B^{out}_{\la_0}$.
% $\wt f[r_{-2\la_0}]$ and $h[r]$, $r_{-2\la_0}$ in
%  $R_{-2\la_0}$ and $r$ in $R$.
\end{prop}

{\em Proof.} For $\rho$ in $R_{(a)}$, set $\wt h[\rho] =
h[(q^\pa\cA)^{-1}\rho]$. Extend $\wt h$ to $\cK$ by linearity.  A
system of relations for $U_\hbar\B_-$ is then
$$
[\wt f[\eps],\wt f[\eps']] = 0, \quad [\wt h[\eps], \wt f[\eta]] = -2 
\wt f[\eps\eta], 
$$
$$
[\wt h[\eps], \wt h[\eps']] = f(\eps,\eps'), 
$$ with $f(\eps,\eps')$ scalar, for $\eps,\eps',\eta$ in $\cK$.
Denote by $\CC\langle \phi, \phi\in F\rangle$ the subalgebra of
$U_\hbar\B_-$ generated by the family $F$ of elements of
$U_\hbar\B_-$.  The product map from
$$ \CC \langle \wt h[o],\wt f[\la'],o\in\mm,\la'\in \La' \rangle \otimes
\CC[\wt h[r_a] ]\otimes \CC \langle \wt h[r],\wt f[r_{-2\la_0}],r\in R,
r_{-2\la_0}\in R_{-2\la_0} \rangle
$$ to $U_\hbar \B_-$ then defines an isomorphism.  Therefore $U_\hbar
\B_-$ is the direct sum of $\CC[\wt h[r_a]]U_\hbar \B^{out}_{\la_0}$
and the left ideal $I$ generated by the $\wt h[r], \wt
f[r_{-2\la_0}]$, $r$ in $R$, $r_{-2\la_{0}}$ in $R_{-2\la_0}$.
$\CC[\wt h[r_a]]U_\hbar \B^{out}_{\la_0}$ is equal to $\CC[h[r_a]]
U_\hbar \B^{out}_{\la_0}$; on the other hand, $\wt f[r_{-2\la_0}] =
\sum_i \res_{P_i}(f(z) k^-(q^{-\pa}z) r_{-2\la_0}(z)\omega_z)$, and
$k^-(q^{-\pa}z)$ belongs to $U_{\hbar,\omega}\G\otimes R_z$, so that
since $R_{-2\la_0}$ is a $R$-module, $\wt f[r_{-2\la_{0}}]$ belongs to
$f[r_{-2\la_{0}}] + \hbar$(right ideal generated by the
$f[\rho_{-2\la_{0}}],\rho_{-2\la_{0}}$ in $R_{-2\la_0}$). Therefore,
$I$ coincides with the left ideal generated by the $h[r],
f[r_{-2\la_0}]$, $r$ in $R$, $r_{-2\la_{0}}$ in $R_{-2\la_0}$, which
is the augmentation ideal of $U_\hbar \B^{out}_{\la_0}$. The Lemma
follows.  \hfill \qed \medskip

\section{A presentation of $U_{\hbar,\omega}\G / (K+2)$} 
\label{other:pres}

Set
$$ \wt e(z) = k^+(q^\pa z)^{-1} k_R(z)^{-1} e(q^\pa z), 
$$ an
\begin{equation} \label{31a} 
  k^{+}_{tot}(z) = k^+(q^{2\pa}z) k_R(q^\pa z)k_R(z)^{-1} k^-(z),
\end{equation}
\begin{equation} \label{31b} 
  k^{-}_{tot}(z) = k^+(q^\pa z)^{-1} k_R(z)^{-1} k_R(q^{-\pa}z)
  k^-(q^{-\pa}z)^{-1}.
\end{equation}

\begin{prop} 
  The following relations
\begin{equation} \label{31c}
  [\wt e(z),\wt e(w)] = [\wt f(z), \wt f(w)] = 0,
\end{equation}
\begin{equation} \label{wte:wtf}
  [\wt e(z), \wt f(w)] = {1\over \hbar} \delta(q^\pa z, w)
  k^+_{tot}(z) - {1\over \hbar} \delta(z,q^\pa w) k^-_{tot}(z)
  {{\exp(2\al(q^{-\pa} z, q^{-\pa}z))} \over {\exp(2\al(q^\pa z,
      q^{-\pa}z))}}
% (K^-(w),k_R(w))_{w = q^{-\pa}z},
\end{equation} 
are satisfied in $U_{\hbar,\omega}\G$.
\end{prop}

{\em Proof.} The proof of $[\wt e(z), \wt e(w)] =0$ is similar to that
of Prop.\ \ref{walsingham}, {\em i)} and relies on the identities
$(k_R(z),e(q^\pa w)) = \exp(2\al(q^\pa w,z))$, and
\begin{equation} \label{identity:alpha:j}
  {{\exp(2\al(q^\pa z,w))}\over {\exp(2\al(q^\pa w,z))}} j(q^\pa z,
  q^\pa w) = 1,
\end{equation} which follows from (\ref{def:al}) and (\ref{id:j}).

Then we have
\begin{align*} 
  \wt f(w) \wt e(z) = & (k_R(w)k^-(q^{-\pa}w), k^+(q^\pa
  z)^{-1}k_R(z)^{-1}) \\ & (f(w), k^+(q^\pa z)^{-1}k_R(z)^{-1})
  (k_R(w)k^-(q^{-\pa}w), e(z)) \\ & k^+(q^\pa z)^{-1}k_R(z)^{-1} f(w)
  e(q^\pa z) k_R(w)k^-(q^{-\pa}w); 
\end{align*}
this equation may be written as 
$$ \wt f_n \wt e_m = \sum_{p\geq 0} \hbar^p \sum_{i\geq N(p), j\geq
  M(p)} A_{ij}^{(p)} k^+_{-i} f_{n-j}e_{m+i} k^-_j, 
$$ where we set $x(z) = \sum_n x_n z^{-n}$, $x = e,f, \wt e, \wt f,
k^{\pm}$; the right side belongs to the completion
$U_{\hbar,\omega}\G$.  Equation (\ref{wte:wtf}) then follows from the
identity
$$ (K^-(q^{-\pa}z),k_R(q^{-\pa}z))= {{\exp(2\al(q^{-\pa} z,
    q^{-\pa}z))} \over {\exp(2\al(q^\pa z, q^{-\pa}z))}}.
$$ \hfill \qed \medskip

\begin{thm}
  $U_{hbar,\omega}\G/(K+2)$ has a presentation with generating series
  $\wt e(z), \wt f(z), k^\pm(z)$ and relations (\ref{17b}),
  (\ref{17c}), (\ref{31a}), (\ref{31b}), (\ref{31c}), (\ref{wte:wtf}),
  and
\begin{equation} \label{k+:wte} 
  ( k^+(z) , \wt e(w)) = q_+(z,q^\pa w) ,
\end{equation}
\begin{equation} \label{k-:wte} 
  ( k^-(z) , \wt e(w)) = \exp[2\al(q^\pa z,w)- 2\al(q^{2\pa} z,w)]
  {{q_+(q^\pa w, q^\pa z)}\over{q_+(q^\pa w,
      q^{2\pa}z)q_-(q^{3\pa}z,q^\pa w)}}
\end{equation} 
and
\begin{equation} \label{k:wtf} 
  ( k^+(z) , \wt f(w)) = q_+(z,q^\pa w)^{-1}, \quad ( k^-(z) , \wt
  f(w)) = ( k^-(z) , \wt e(w))^{-1}.
\end{equation} 
\end{thm}

\section{Central current $T(z)$}

Recall that
$$
g^+_{\la}(z) = (G_{-2\la} - G)(q^\pa z,z), \quad 
g^-_{\la}(z) = (G_{-2\la} - G)(q^{-\pa}z,z). 
$$

Define $\sigma,\al,\beta$ in $R[[\hbar]]$ and $A_\la,B_\la$ in
$\cK[[\hbar]]$ by
\begin{equation} \label{sigma}
 \sigma(q^\pa z) = \left[ - e^{-2\sum_i (q^{\pa}U_+ e_i)(z) \otimes
  e^i(w)} e^{-\phi(\hbar,\pa_z^i \gamma)} \psi(\hbar,\pa_z^i\gamma)
\right]_{z=w}; 
\end{equation} 
\begin{equation} \label{alpha}
\al(q^\pa z) = \left[ - e^{- 2\sum_i (q^\pa U_+ e_i)(z) \otimes e^i(w)}
\pa_{\hbar} \{e^{-\phi(\hbar,\pa_z^i\gamma)}\psi(\hbar,\pa_z^i\gamma)\}
\right]_{w = z} , 
\end{equation}
\begin{equation} \label{beta}
  \beta(q^\pa z) = \al(q^\pa z) - 2 \pa_\hbar[\tau_{w=z}] \si(q^\pa z), 
\end{equation}
\begin{equation} \label{A:la}
  A_\la(z) = \al(q^{2\pa}z) + \sigma(q^{2\pa}z) [g^+_{\la}(z) - \sum_i
  e^i (z) (q^{2\pa}(q^{-\pa}e_i)_R)(z) ] , 
\end{equation}
\begin{equation} \label{B:la}
  B_\la(z) = \beta(q^{2\pa}z) - \sigma(q^{2\pa}z) [g^-_{\la}(z) -
  \sum_i e^i(z) ((q^{-\pa}e_i)_R)(z)] ,
\end{equation}
we have $\sigma(z) = \hbar + O(\hbar^2)$, $\al(z) = 1 + O(\hbar)$,
$\beta(z) = 1 + O(\hbar)$, $A_\la = 1+ O(\hbar)$, $B_\la = 1+
O(\hbar)$.

Let us set 
$$ T(z) = \wt e(z) \wt f(z)_{z\to R_{2\la}} + \wt f(z)_{z\to \La'}\wt
e(z) + a_\la(z) k^{+}_{tot}(z) + b_\la(z) k^{-}_{tot}(z) , 
$$
where
\begin{equation} \label{a:la}
  a_\la(z) = {1\over \hbar} {{A_\la(z)}\over{\sigma(q^{2\pa}z)}},
\end{equation}
and
\begin{equation} \label{b:la}
b_\la(z) = {1\over
  \hbar}{{\exp(2\al(q^{-\pa}z,q^{-\pa}z))}\over{\exp(2\al(q^\pa z,
    q^{-\pa}z))}} {{B_\la(z)}\over{\sigma(q^{2\pa}z)}};
\end{equation}
we will also set
\begin{equation} \label{b':la}
  b'_\la(z) = {1\over \hbar}{{B_\la(z)}\over{\sigma(q^{2\pa}z)}}.
\end{equation}

\begin{thm} \label{Sugawara}
  The Laurent coefficients of $T(z)$ are central elements of
  $U_{\hbar,\omega}\G$.
\end{thm}

The proof is contained in the next sections.

\subsection{Commutation of $T(z)$ with $\wt e(w)$}

Set
$$ :k^+_{tot}(z) \wt e(w): = k^+(q^{2\pa}z) k^+(q^\pa w)^{-1} k_R(q^\pa z)
k_R(z)^{-1} k_R(w)^{-1} \wt e(q^\pa w) k^-(z) ,$$ and
$$ :k^-_{tot}(z) \wt e(w): = k^+(q^\pa z)^{-1} k_R(z)^{-1}
k_R(q^{-\pa}z) k^+(q^\pa w)^{-1} k_R(w)^{-1} \wt e(q^\pa w)
k^-(q^{-\pa}z)^{-1} . $$

\begin{lemma} \label{pdt:k+:wte}
  We have
$$ k^+_{tot}(z) \wt e(w) = \exp(2\al(q^\pa w,z)- 2\al(q^\pa w,q^\pa
z)) q_-(q^{2\pa}z,q^\pa w)^{-1} :k^+_{tot}(z) \wt e(w): , 
$$ and
$$ \wt e(w) k^+_{tot}(z) = \exp(2\al(q^\pa w,z)- 2\al(q^\pa w,q^\pa
  z)) q_+(q^{2\pa}z,q^\pa w)^{-1} :k^+_{tot}(z) \wt e(w): . 
  $$
\end{lemma}

{\em Proof.} Let us prove the first identity. The factor in the right
side is
\begin{equation} \label{first:factor}
(k^-(z),k^+(q^\pa w)^{-1}) (k^-(z),k_R(w)^{-1})(k^-(z),e(q^\pa w)). 
\end{equation}
we have 
$$ (k^-(z),k_R(w)^{-1}) = \exp(2\al(q^\pa z,w)- 2\al(q^{2\pa} z,w)),
$$
therefore (\ref{first:factor}) is equal to 
\begin{equation} \label{panthoos}
{{q_+(q^\pa w, q^\pa z)}\over{q_+(q^\pa w, q^{2\pa}z)}}
\exp(2\al(q^\pa z, w)- 2\al(q^{2\pa} z, w) ) q_-(q^{3\pa}z,q^\pa
w)^{-1}. 
\end{equation}
Identity (\ref{identity:alpha:j}) can be formulated as 
$$
{{\exp(2\al(q^\pa w,z))}\over{\exp(2\al(q^\pa z,w))}} = q_+(q^\pa w,q^\pa z)
q_-(q^{2\pa}z,q^\pa w), 
$$ because the right side is $j(q^\pa z, q^\pa w)$ (see (\ref{j})).
Applying to this identity $\exp\circ (1 - q^{\pa_z}) \circ \log$, 
we transform (\ref{panthoos}) into 

$$ \exp(2\al(q^\pa w,z)- 2\al(q^\pa w,q^\pa z))
q_-(q^{2\pa}z,q^\pa w)^{-1};
$$
this implies the first equality. 

The factor in the right side of the second identity is
$$
(e(q^{\pa}w), k^+(q^{2\pa}z)k_R(q^\pa z)k_R(z)^{-1}), 
$$
which is equal to 
$$ \exp(2\al(q^\pa w,z)- 2\al(q^\pa w,q^\pa z)) q_+(q^{2\pa}z,q^\pa
w)^{-1}.
$$ \hfill \qed \medskip

In the same way, one proves 
\begin{lemma} \label{pdt:k-:wte}
We have 
$$ k^-_{tot}(z) \wt e(w) = \exp(2\al(q^\pa w, z) -2\al(q^\pa w,
q^{-\pa} z)) q_-(q^\pa z, q^\pa w) :k^-_{tot}(z) \wt e(w): ,
$$ and
$$ \wt e(w) k^-_{tot}(z) = \exp(2\al(q^\pa w, z) -2\al(q^\pa w,
q^{-\pa} z)) q_+(q^\pa z, q^\pa w) :k^-_{tot}(z) \wt e(w): .
$$
\end{lemma}

Then we have 

\begin{prop} \label{comm:e}
$T(z)$ commutes with $\wt e(w)$. 
\end{prop}

{\em Proof.} We have 
\begin{align} \label{[T,wte]}
  & \nonumber [T(z),\wt e(w)] = \wt e(z) [\wt f(z), \wt e(w)]_{z\to
    R_{2\la}} + [\wt f(z), \wt e(w)]_{z\to \La'} \wt e(z) \\ &
  \nonumber + a_\la(z) [k^+_{tot}(z), \wt e(w)] + b_\la(z)
  [k^-_{tot}(z), \wt e(w)] \\ & \nonumber = - \wt e(z) \left(
    {1\over\hbar} \delta(q^\pa w,z) k^+_{tot}(w) - {1\over \hbar}
    \delta(w,q^\pa z) k^-_{tot}(w)
    {{\exp(2\al(q^{-\pa}w,q^{-\pa}w))}\over{\exp(2\al(q^\pa w,
        q^{-\pa}w))}} \right)_{z\to R_{2\la}} \\ & \nonumber - \left(
    {1\over\hbar} \delta(q^\pa w,z) k^+_{tot}(w) - {1\over \hbar}
    \delta(w,q^\pa z) k^-_{tot}(w)
    {{\exp(2\al(q^{-\pa}w,q^{-\pa}w))}\over{\exp(2\al(q^\pa w,
        q^{-\pa}w))}} \right)_{z\to \La'} \wt e(z) \\ & \nonumber +
  a_\la(z) [k^+_{tot}(z), \wt e(w)] + b_\la(z) [k^-_{tot}(z), \wt
  e(w)] \\ & \nonumber = - {1\over \hbar} \left( G_{-2\la}(q^\pa w,z)
    k^+_{tot}(w)\wt e(z) + G_{2\la}(z,q^\pa w) \wt e(z) k^+_{tot}(w)
  \right) \\ & \nonumber + {1\over\hbar} \left( G_{-2\la}(q^{-\pa}
    w,z) k^-_{tot}(w) \wt e(z) + G_{2\la}(z,q^{-\pa} w) \wt e(z)
    k^-_{tot}(w) \right)
  {{\exp(2\al(q^{-\pa}w,q^{-\pa}w))}\over{\exp(2\al(q^\pa w,
      q^{-\pa}w))}}\\ & + a_\la(z) [k^+_{tot}(z), \wt e(w)] + b_\la(z)
  [k^-_{tot}(z), \wt e(w)];
\end{align}
the last equality follows from the identities 
$$ \delta(w,z)_{z\to R_{2\la}} = G_{2\la}(z,w), \quad
\delta(w,z)_{z\to \La'} = G_{-2\la}(w,z).
$$

We have 
\begin{align*}
& G_{-2\la}(q^\pa w,z) k^+_{tot}(w)\wt e(z) + G_{2\la}(z,q^\pa w) \wt
e(z) k^+_{tot}(w)
\\ & = 
\exp( 2\al(q^\pa z, w)- 2\al(q^\pa z, q^\pa w)) \\ & \left(
    G_{-2\la}(q^\pa w,z) q_-(q^{2\pa}w,q^\pa z)^{-1} +
    G_{2\la}(z,q^\pa w) q_+(q^{2\pa}w,q^\pa z)^{-1} \right) 
:k^+_{tot}(w)\wt e(z):
\\ & = 
 \exp( 2\al(q^\pa z, w)- 2\al(q^\pa z, q^\pa w)) A_\la(z)
  \delta(z,w):k^+_{tot}(w)\wt e(z):,  
\end{align*}
where the first equality follows from Lemma \ref{pdt:k+:wte}, 
and the second from Lemma \ref{ids:G:q}. 

In the same way, we have 
\begin{align*}
  & G_{-2\la}(q^{-\pa} w,z) k^-_{tot}(w) \wt e(z) +
  G_{2\la}(z,q^{-\pa} w) \wt e(z) k^-_{tot}(w) \\ & = 
  \exp(2\al(q^\pa z, w) -2\al(q^\pa z, q^{-\pa} w)) \\ & \left(
    G_{-2\la}(q^{-\pa} w,z) q_-(q^\pa w, q^\pa z) +
    G_{2\la}(z,q^{-\pa} w) q_+(q^\pa w, q^\pa z) \right) :k^-_{tot}(z)
  \wt e(w) :\\ & = \exp(2\al(q^\pa z, w) -2\al(q^\pa z, q^{-\pa} w)) B_\la(z)
  \delta(z,w):k^-_{tot}(z)
  \wt e(w) :
\end{align*}
where the first equality follows from Lemma \ref{pdt:k-:wte}, and the
second from Lemma \ref{ids:G:q}.

%$$G_{-2\la}(q^\pa w,z) k^+_{tot}(w)\wt e(z) + G_{2\la}(z,q^\pa w) \wt
%e(z) k^+_{tot}(w)
%$$
%is equal to 
%\begin{align*}
%  & \exp( 2\al(q^\pa z, w)- 2\al(q^\pa z, q^\pa w)) \\ & \left(
%    G_{-2\la}(q^\pa w,z) q_-(q^{2\pa}w,q^\pa z)^{-1} +
%    G_{2\la}(z,q^\pa w) q_+(q^{2\pa}w,q^\pa z)^{-1} \right) \\ &
%  k^+(q^{2\pa}w)k^+(q^\pa z)^{-1} k_R(q^\pa w) k_R(w)^{-1} k_R(z)^{-1}
%  e(q^\pa z) k^-(z); 
%\end{align*}
%by Lemma \ref{ids:G:q}, this is equal to
%\begin{align*}
%  & \exp( 2\al(q^\pa z, w)- 2\al(q^\pa z, q^\pa w)) A_\la(z)
%  \delta(z,w) \\ & k^+(q^{2\pa}w)k^+(q^\pa z)^{-1} k_R(q^\pa w)
%  k_R(w)^{-1} k_R(z)^{-1} e(q^\pa z) k^-(z)
%\end{align*}
%and by Lemma \ref{pdt:k-:wte},
%$$ G_{-2\la}(q^{-\pa} w,z) k^-_{tot}(w) \wt e(z) + G_{2\la}(z,q^{-\pa}
%w) \wt e(z) k^-_{tot}(w)
%$$
%is equal to 
%\begin{align*}
%  & \exp(2\al(q^\pa z, w) -2\al(q^\pa z, q^{-\pa} w)) \\ & \left(
%    G_{-2\la}(q^{-\pa} w,z) q_-(q^\pa w, q^\pa z) +
%    G_{2\la}(z,q^{-\pa} w) q_+(q^\pa w, q^\pa z) \right) \\ &
%  k^+(q^\pa w)^{-1} k_R(w)^{-1} k_R(q^{-\pa}w) k^+(q^\pa z)^{-1}
%  k_R(z)^{-1} e(q^\pa z) k^-(q^{-\pa}w)^{-1}; 
%\end{align*}
%by Lemma \ref{ids:G:q}, this is is equal to
%\begin{align*}
%  & \exp(2\al(q^\pa z, w) -2\al(q^\pa z, q^{-\pa} w)) B_\la(z)
%  \delta(z,w) \\ & k^+(q^\pa w)^{-1} k_R(w)^{-1} k_R(q^{-\pa}w)
%  k^+(q^\pa z)^{-1} k_R(z)^{-1} e(q^\pa z) k^-(q^{-\pa}w)^{-1} .
%\end{align*}

On the other hand, we have
\begin{align*}
  & [k^+_{tot}(z), \wt e(w)] \\ & = \exp(2\al(q^\pa w,z)- 2\al(q^\pa
  w,q^\pa z)) [q_-(q^{2\pa}z,q^\pa w)^{-1} - q_+(q^{2\pa}z,q^\pa
  w)^{-1}] \\ & 
:k^+_{tot}(z) \wt e(w):
\\ & = \exp(2\al(q^\pa
  w,z)- 2\al(q^\pa w,q^\pa z)) \sigma(q^{2\pa}z) \delta(z,w)
:k^+_{tot}(z) \wt e(w):
\end{align*}
and
\begin{align*}
  & [k^-_{tot}(z) , \wt e(w)] \\ & = \exp(2\al(q^\pa w,z) - 2\al(q^\pa
  w, q^{-\pa}z)) [q_-(q^\pa w, q^\pa z)- q_+(q^\pa w, q^\pa z) ] \\ &
  :k^-_{tot}(z) \wt e(w): \\ & = \exp(2\al(q^\pa w,z) - 2\al(q^\pa w,
  q^{-\pa}z)) [-\sigma(q^{2\pa}w)\delta(z,w)] :k^-_{tot}(z) \wt e(w):
  .
\end{align*}

The equalities $\delta(z,w) :k^{\pm}_{tot}(z) \wt e(w):= \delta(z,w)
:k^{\pm}_{tot}(w) \wt e(z):$ then imply that 
(\ref{[T,wte]}) vanishes. \hfill \qed

\subsection{Commutation of $T(z)$ with $k^\pm(w)$}

Let us denote by $U_\hbar\N_+$, $U_\hbar\HH$ and $U_\hbar\N_-$ the
subalgebras of $U_{\hbar,\omega}\G$ generated respectively by the
$e[\eps]$, by the $h^+[r]$, $h^-[\la]$ and $K$, and by the $f[\eps]$.
If we assign degree $1$ to the $e[\eps]$ and $f[\eps]$,
$U_\hbar\N_\pm$ are graded algebras. We denote by
$U_\hbar\N_\pm^{[i]}$ their homogeneous components of degree $i$.

We will prove 

\begin{lemma} \label{bill}
  $k^+(w) T(z) k^+(w)^{-1} - T(z)$ and $k^-(w) T(z) k^-(w)^{-1} -
  T(z)$ both belong to $U_\hbar\HH$.
\end{lemma}

{\em Proof.} It suffices to prove the same statements with $T(z)$
replaced by $T_0(z)$ defined by
$$
T_0(z) =  \wt e(z) \wt f(z)_{z\to R_{2\la}} + \wt f(z)_{z\to \La'}\wt
e(z) .
$$ Then from (\ref{k+:wte}) and (\ref{k:wtf}) follows that
\begin{align*}  
  & q_+(z,q^\pa w)^{-1} k^+(w) T_0(z) k^+(w)^{-1} - T_0(z) \\ & = \wt
  e(w) [q_+(z,q^\pa w)^{-1} \wt f(w)]_{w\to R_{2\la}} + [q_+(z,q^\pa
  w)^{-1} \wt f(w)]_{w\to \La'} \wt e(w) \\ & - q_+(z,q^\pa w)^{-1}
  [\wt e(w) \wt f(w)_{w\to R_{2\la}}+ \wt f(w)_{w\to \La'} \wt e(w) ]
  \\ & = [\wt e(w) , [q_+(z,q^\pa w)^{-1} \wt f(w)]_{w\to \La'} -
  q_+(z,q^\pa w)^{-1} \wt f(w)_{w\to \La'} ]
\end{align*}
because of the identity 
\begin{align*} 
  & [q_+(z,q^\pa w)^{-1} \wt f(w)]_{w\to \La'} - q_+(z,q^\pa w)^{-1}
  \wt f(w)_{w\to \La'} \\ & = q_+(z,q^\pa w)^{-1} \wt f(w)_{w\to
    R_{2\la}} - [q_+(z,q^\pa w)^{-1} \wt f(w)]_{w\to R_{2\la}} .
\end{align*} 
For any $\eps$ in $\cK$, $[\wt e[\eps],\wt f(z)]$ belongs to
$U_\hbar\HH$, which proves the first part of the statement. The second
part is proved in the same way, using (\ref{k+:wte}) and
(\ref{k:wtf}).  \hfill \qed \medskip

Let us now prove 
\begin{prop} \label{comm:h}
$T(z)$ commutes with $U_\hbar\HH$. 
\end{prop}

{\em Proof.} Set for $r$ in $R$ and $\la$ in $\La$, 
$$ x^+_\eta(r) = [h^+[r],T[\eta]], \quad x^-_\eta(\la) =
[h^-[\la],T[\eta]].
$$ From Lemma \ref{bill} follows that $x^\pm_\eta$ are linear maps
from $R$ and $\La$ to $U_\hbar\HH$. Moreover, we have $[x^+_\eta(r) ,
\wt f[\eps]] = [[h^+[r],T[\eta]],\wt f[\eps]] = - [[T[\eta],\wt
f[\eps]], h^+[r]] - [[\wt f[\eps],h^+[r]],T[\eta]]$; both terms are zero
by Prop.\ \ref{comm:f}, so that we have
$$
[x^+_\eta(r) , \wt f[\eps]] = 0; 
$$ in the same way, one shows that
$$
[x^-_\eta(\la) , \wt f[\eps]] = 0. 
$$ But any element $x$ of $U_\hbar\HH$, such that $[x,\wt f[\eps]] =
0$ for any $\eps$, is zero. To show this, one may divide $x$ by the
greatest possible power of $\hbar$ and check that the same statement
is true in the classical affine Kac-Moody algebra. \hfill \qed

\subsection{Commutation of $T(z)$ with $f(w)$}

\begin{lemma}
  $T(z)$ may be written
  $$ T(z) = \wt f(z) \wt e(z)_{z\to R} + \wt e(z)_{z\to \La} \wt f(z)
  + \kappa(z),
  $$ where $\kappa(z)$ belongs to $U_\hbar\HH[[z,z^{-1}]]$.
\end{lemma}

{\em Proof.} We have
$$ \wt f(z) \wt e(z)_{z\to R} + \wt e(z)_{z\to \La} \wt f(z) - T_0(z)
$$ 
is equal to 
$$ [\wt e(z)_{z\to \La'}, \wt f(z)_{z\to \La'}] - [\wt e(z)_{z\to
  R_{2\la}}, \wt f(z)_{z\to R_{2\la}}]
$$ and therefore belongs to $U_\hbar\HH[[z,z^{-1}]]$.  \hfill \qed
\medskip

We first show: 

\begin{lemma}
  The commutator $[T(z),\wt f(w)]$ belongs to $U_\hbar \HH U_\hbar
  \N_+^{[1]}$; in other words, there are formal series $K_i(z,w)$ in
  $U_\hbar\HH[[z,z^{-1},w,w^{-1}]]$, such that
  \begin{equation} \label{balus}
    [T(z),\wt f(w)] = \sum_i K_i(z,w) \wt f[\eps_i].
  \end{equation}
\end{lemma}

{\em Proof.} It suffices to show this with $T_0(z)$ instead of $T(z)$.
This follows from a reasoning analogous to the first part of the proof
of Prop.\ \ref{comm:e}. 
\hfill \qed\medskip

{}From there follows: 
\begin{prop} \label{comm:f}
  $T(z)$ commutes with $\wt f(w)$.
\end{prop}

{\em Proof.} Let $\eps$ belong to $\cK$. $\wt e[\eps]$ commutes with the
left side of (\ref{balus}), by Props.\ \ref{comm:f} and \ref{comm:e}.
Let us write that it commutes with the right side of this equality. We
get $\sum_i [\wt e[\eps],K_i(z,w)] \wt f[\eps_i] + $ element of $U_\hbar\HH =
0$. From there follows that $[\wt e[\eps],K_i(z,w)] = 0$. The reasoning of
the end of the proof of Prop.\ \ref{comm:h} applies to show that
$K_i(z,w)$ vanishes. \hfill \qed \medskip

Props.\ \ref{comm:e}, \ref{comm:h} and \ref{comm:f} imply Thm.\
\ref{Sugawara}.  \hfill \qed \medskip

\begin{remark} {\it Classical limit.}
  Let us show that $T(z)$ is, up to a scalar, a deformation of the
  Sugawara tensor.  Let us denote by $e_{cl}(z),h_{cl}(z)$ and
  $f_{cl}(z)$ the generating currents of $\G$. Then we have
  $$ e(z) = e_{cl}(z) + O(\hbar), \quad f(z) = f_{cl}(z) + O(\hbar),
  $$
  $$ k^+(z) = 1 + {\hbar \over 2} h_{cl}(z)_{z\to \La} + o(\hbar),
  \quad k^-(z) = 1 + {\hbar \over 2} h_{cl}(z)_{z\to R} + o(\hbar),
  $$ 
$k_R(z) = 1 + O(\hbar^2)$, so that 
$$ k^+_{tot}(z) = [1 + {\hbar \over 2} q^{2\pa_z}(h_{cl}(z)_{z\to
  \La}) + \hbar^2 s(z)] [1 + {\hbar \over 2} h_{cl}(z)_{z\to R} +
\hbar^2 t(z)] + O(\hbar^3)
$$ 
\begin{align*} & k^-_{tot}(z) =  [1 - {\hbar \over 2}
  q^{\pa_z}(h_{cl}(z)_{z\to \La}) - \hbar^2 s(z) + {\hbar^2\over
    4}(h_{cl}(z)_{z\to \La})^2] \\ & [1 - {\hbar \over 2}
  q^{-\pa_z}(h_{cl}(z)_{z\to \La}) - \hbar^2 t(z) + {\hbar^2\over
    4}(h_{cl}(z)_{z\to R})^2] + O(\hbar^3),
\end{align*} 
where $s(z)$ and $t(z)$ are some currents. Then
\begin{align*}
  & T(z) = e_{cl}(z)_{z\to \La}f_{cl}(z) + f_{cl}(z) e_{cl}(z)_{z\to
    R} + {1\over{\hbar^2}}(k^+_{tot}(z) + k^-_{tot}(z)) + O(\hbar) \\ 
  & = {1\over{\hbar^2}} + e_{cl}(z)_{z\to \La}f_{cl}(z) + f_{cl}(z)
  e_{cl}(z)_{z\to R} + {1\over 2} \pa h_{cl}(z) \\ & + {1\over 4}
  \left( h_{cl}(z)_{z\to \La} h_{cl}(z)+ h_{cl}(z) h_{cl}(z)_{z\to R}
  \right) + O(\hbar);
\end{align*}
so $T(z) - \hbar^{-2}$ coincides with the classical Sugawara tensor to
order $\hbar$. \hfill \qed \medskip
\end{remark}

\begin{remark} {\it Other expressions of $T(z)$.}\label{old:T}
  One may show that up to an additive scalar constant, $T(z)$
  coincides with
\begin{align} \label{Sug:la} 
  T'(z) = & k^+(q^\pa z)^{-1} \left( f(z)_{z\to \La'} e(q^\pa z) +
    e(q^\pa z)f(z)_{z\to R_{2\la}} \right) k^-(q^{-\pa}z) \\ &
  \nonumber + {{\gamma'_{\la}(z)}\over{\hbar\sigma(z)}} k^+(q^\pa
  z)^{-1} k^-(q^{-2\pa}z)^{-1} +
  {{\delta_{2\la}(z)}\over{\hbar\sigma(q^{2\pa} z) }} k^+(q^{2\pa} z)
  k^-(q^{-\pa}z) ,
\end{align}
with 
$$ \gamma'_\la(z) = \al(q^{2\pa}z) - \sigma(q^{2\pa}z)[ g^-_\la(z) +
\sum_i e^i(z) (q^\pa (q^{-2\pa}e_i)_R)(z) ]
$$
and 
$$
\delta_\la(z) = \beta(q^{2\pa}z) + \sigma(q^{2\pa}z) g^+_\la(z). 
$$ 
It also coincides with $T''(z)$ defined by 
\begin{align} \label{T''} 
  T''(z) & = k^+(z)^{-1} \left( e(z)_{z\to \La}f(q^{-\pa}z) +
    f(q^{-\pa}z)e(z)_{z\to R} \right) k^-(q^{-2\pa}z) \\ &
  \nonumber + {1\over\hbar}{\al\over \sigma}(z)
  k^+(q^{-\pa}z)k^-(q^{-2\pa}z) + {1\over \hbar}{\beta' \over
    \sigma}(z) k^+(z)^{-1} k^-(q^{-\pa}z)^{-1},
\end{align} 
up to an additive constant, with $\beta'(z) = \beta(q^\pa z) -
\sigma(q^\pa z) \sum_i q^\pa((q^{-2\pa}e_i)_R)(z) e^i(z)$.

This formula is a generalization of the formula given in \cite{FR},
which uses \cite{RS} and the new realizations isomorphism. To see the
correspondance between this formula and ours, let us modify the
notation in \cite{FR} so that the quantum parameter of that paper is
denoted by $\underline{q}$.  The level in \cite{FR} is denoted by $k$
and the currents generating the algebra $U_q\G$ are $k_1^\pm(z)$,
$E(z)$ and $F(z)$.

Set $X = \CC P^1$ and $\pa = z{d\over{dz}}$. The algebra
$U_{\hbar,\omega}\G$ is isomorphic to $U_q\G$, the isomorphism $i$
being given by the formulas
$$ i(K) = k, \quad i(k^+(z)) = k_1^+(z\underline{q}^{{k\over
    2}+2})^{-1}, \quad i(k^-(z)) = k_1^-(z\underline{q}^{3k\over
  2}),$$
$$ i(e(z)) = - {1\over{\hbar(\underline{q}-\underline{q}^{-1})}} E(z),
i(f(z)) = F(\underline{q}^k z),
$$ with
$$ q = \underline{q}^{-2}, \quad q(z,w) =
{{q^{-1}z-w}\over{z-q^{-1}w}}, \quad q_+(z,w) =
{{q^{-1/2}z-q^{1/2}w}\over{z-w}}.
$$

Formula (6.10) of \cite{FR} then gives
\begin{align*}
  & i^{-1} ( \ell(z) ) =
  {1\over{\hbar(\underline{q}-\underline{q}^{-1})}} k^+(z)^{-1}
  :e(z)f(q^{-1}z): k^-(zq^{-2}) \\ & + q^{-1/2}
  k^+(zq^{-1})k^-(zq^{-2})+ q^{1/2} k^+(z)^{-1}k^-(zq^{-1})^{-1},
\end{align*}
so $i^{-1}(\ell(z))$ is equal to $T(z)$ given by (\ref{T''}).  \hfill
\qed \medskip
\end{remark}

\begin{remark} \label{genus:one}
 {\it Genus $1$ case.}
  Assume $X$ is an elliptic curve $\CC/L$, $L = \ZZ + \tau\ZZ$, and
  $\omega = dz$. Let $\theta$ be the Jacobi theta-function, equal to 
  $$ \theta(z) = {{\on{sin}(\pi z)}\over{\pi}} \prod_{j=1}^{\infty}
  {{(1 - e^{2i\pi(j\tau + z)})(1 - e^{2i\pi(j\tau - z)}) }\over{(1 -
      e^{2i\pi j \tau})^2}}
  $$ The Weierstrass function is $\wp = -(d/dz)^2 \ln \theta(z)$.
  According to Prop.\ \ref{Lambda}, we have $R =\CC 1 \oplus
  (\oplus_{i\geq 0}\CC (d/dz)^i \wp)$ and $\La = \CC\theta'/\theta
  \oplus z\CC[[z]]$. We have also $R_\la = \oplus_{i\geq 0}
  \CC (d/dz)^i({{\theta(z-\la)}\over{\theta(z)}})$ and $\La' = \CC[[z]]$.
  We have
  $$ G(z,w) = d/dz\ln\theta(z-w) - d/dz\ln\theta(z) + d/dz\ln\theta(w) , 
  $$
  $$ G_{2\la}(z,w) = {{\theta(-2\la +
      z-w)}\over{\theta(z-w)\theta(-2\la)}}, $$
%  $$ so that
%  $$ \tilde g_{2\la}(z) = - {\theta'\over\theta}(2\la).
%  $$ We have also 
  $q_-(z,w) = {{\theta(z-w-\hbar)}\over{\theta(z-w)}}$, viewed as a
  series in $\CC((z))((w))[[\hbar]]$,
  $$ \sigma(z) = \theta(\hbar),
%  $$ \al(z) = \theta'(\hbar) +
%  \theta(\hbar)\left({\theta'\over\theta}(z-\hbar) -
%    {\theta'\over\theta}(z) \right), \beta(z) = \theta'(\hbar)
%  -\theta(\hbar) \left( {\theta'\over\theta}(z-\hbar) -
%    {\theta'\over\theta}(z-2\hbar) \right),  $$ 
  \quad \gamma'_{\la}(z) = {{\theta(2\la-\hbar)}\over{\theta(2\la)}},
  \quad \delta_{\la}(z) = {{\theta(2\la+\hbar)}\over{\theta(2\la)}}.
$$
The expression of $T'(z)$ is then 
\begin{align*}
  T'(z) = & k^+(z+\hbar)^{-1} (f(z)_{z\to\La'} e(z+\hbar) + e(z+\hbar)
  f(z)_{z\to R_{2\la}}) k^-(z-\hbar) \\ & + {{\theta(2\la -
      \hbar)}\over{\hbar\theta(\hbar)\theta(2\la)}}
k^+(z+\hbar)^{-1} k^-(z-2\hbar)^{-1} %\\ & 
+{{\theta(2\la+\hbar)}\over{\hbar \theta(2\la)\theta(\hbar)}}
k^+(z+2\hbar) k^-(z-\hbar); 
\end{align*}
we have also 
$$ T(z) = \wt e(z) \wt f(z)_{z\to R_{2\la}}+ \wt f(z)_{z\to \La'} \wt
e(z) + {1\over\hbar}{{\theta(2\la +
    \hbar)}\over{\theta(2\la)\theta(\hbar)}} k^+_{tot}(z) +
{1\over\hbar}{{\theta(2\la - \hbar)}\over{\theta(2\la)\theta(\hbar)}}
k^-_{tot}(z) .
$$

\hfill \qed
\end{remark}

\section{Subalgebras $U_\hbar \G^{out}$ and $U_\hbar \G^{out}_{\la_0}$ of
  $U_{\hbar,\omega}\G$ and coproducts} \label{subalg}

\subsection{Subalgebras $U_\hbar\G^{out}$ and $U_\hbar\G^{out}_{\la_0}$}

In \cite{ER:qH}, we showed that $U_{\hbar,\omega}\G$ contains a
``regular'' subalgebra $U_\hbar\G^{out}$, generated by the
$h^+[r],e[r]$ and $f[r]$, for $r$ in $R$. The inclusion
$U_\hbar\G^{out} \subset U_{\hbar,\omega}\G$ is a deformation of the
inclusion of the classical enveloping algebra of $\bar\G \otimes R$ in
that of $\G = (\bar\G \otimes \cK) \oplus \CC K$.

For any $\la_0$ in $\CC^g$, define $U_\hbar\G^{out}_{\la_0}$ as the
subalgebra of $U_{\hbar,\omega}\G^{out}$, generated by the $h^+[r]$,
$e[r_{2\la_0}]$ and $f[r_{-2\la_0}]$, for $r$ in $R$, $r_{\pm 2\la_0}$
in $R_{\pm 2\la_0}$.

\begin{prop}
  Define $\G^{out}_{\la_0}$ to be the Lie algebra $(\bar\N_+ \otimes
  R_{2\la_0}) \oplus (\bar\HH \otimes R) \oplus (\bar \N_- \otimes
  R_{-2\la_0})$. The inclusion $U_\hbar\G^{out}_{\la_0} \subset
  U_{\hbar,\omega}\G$ is a deformation of the inclusion of the
  classical enveloping algebra of $\G^{out}_{\la_0}$ in that of $\G$.
\end{prop}

{\em Proof.} (\ref{e-e}) implies that the $e[\eps]$ satisfy the
relations given by the pairing of
$$ (1+\psi( - \hbar, \pa^i_z \gamma)G(z,w)) e(z) e(w) = e^{2(\tau -
  \phi)}(1+\psi(\hbar, \pa^i_z\gamma)G(z,w)) e(w) e(z) $$ with any
$\upsilon$ in $\cK \otimes \cK$, such that $m(\upsilon)=0$, where $m$
is the multiplication map.

Taking for $\upsilon$ any $\al\otimes\beta - \beta \otimes\al$, with
$\al,\beta$ in $R_{2\la_0}$, and using the fact that $R_{2\la_0}$ is
an $R$-module, we get relations of the form 
$$ [e[\al],e[\beta]] = \sum_{i\geq 1,j} \hbar^i
e[\al^{(i)}_j]e[\beta^{(i)}_j] , 
$$ with $\al^{(i)}_j,\beta^{(i)}_j$ in $R_{2\la_0}$. Therefore, if
$e_{2\la_0;i}$ is a basis of $R_{2\la_0}$, the family
$$(e[e_{2\la_0;i_1}]\cdots e[e_{2\la_0;i_p}])_{i_1\leq \cdots \leq
  i_p}$$ spans the subalgebra of $U_{\hbar,\omega}\G$ generated by the
$e[r]$, $r$ in $R_{2\la_0}$.  Since by \cite{ER:qH}, Lemma 3.3, this
is also a free family, it forms a basis of this subalgebra. To finish
the proof, one proves the similar basis result for the subalgebra
generated by the $f[r]$, $r$ in $R_{-2\la_0}$ and a triangular
decomposition result (see \cite{ER:qH}, Prop.\ 3.2 and Prop.\ 3.5).
\hfill \qed \medskip

\subsection{Coproducts} 

Set $A = U_{\hbar,\omega}\G$, $B = U_\hbar \G^{out}$, $B_{\la_0} =
U_\hbar \G^{out}_{\la_0}$.

Define for $\nn = (n_i)_{1\leq i \leq p}$, $I_\nn$ as the left ideal
of $A$ generated by the $x[\eps],\eps\in \prod_i z_i^{n_i}
\CC[[z_i]]$.  Define $A\otimes_> A$, $A\otimes_< A$ and $A \bar\otimes
A$ as the completions of $A\otimes A$ with respect to the topologies
defined by $A\otimes I_\nn$, $I_\nn \otimes A$ and $I_\nn \otimes A
+ A \otimes I_\nn$ ($\otimes$ denotes the $\hbar$-adically completed
tensor product). We have the inclusions $A\otimes_> A \subset A
\bar\otimes A$, $A\otimes_< A \subset A \bar\otimes A$ and $A\otimes A
= (A\otimes_> A) \cap (A\otimes_< A)$.

We define also for any space $V$, $V\otimes_> A$ as the completion of
$V \otimes A$ w.r.t. the topology defined by the $V \otimes I_\nn$,
$A^{\otimes_> n}$ as $A^{\otimes_> n-1}\otimes_> A$, and $A^{\otimes_<
  n}$ in the same way.

In \cite{ER:qH}, we defined Drinfeld-type coproducts $\Delta$ and
$\bar\Delta$ on $U_{\hbar,\omega}\G$ by formulas similar to those of
\cite{Dr:new}.  $\Delta$ and $\bar\Delta$ map $A$ to $A\otimes_< A$
and to $A\otimes_> A$. Moreover, $\Delta$ and $\bar\Delta$ are
conjugated by an element $F$ of $A \bar\otimes A$. $F$ is decomposed
as a product $F_2 F_1$, with $F_1$ in $A\otimes_< B$ and $F_2$ in
$B\otimes_> A$, which are defined as $\lim_{\leftarrow}A \otimes
B/I_\nn \otimes B$ and $\lim_{\leftarrow}B \otimes A/B \otimes I_\nn$.

$\Delta_R$ is defined as $\Ad(F_1) \circ \Delta$. It maps therefore
$A$ to $A\otimes_< A$. Since $\Delta_R$ is equal to $\Ad(F_2^{-1})
\circ \bar\Delta$, it also maps $A$ to $A \otimes_> A$ and therefore to
$A \otimes A$. Also we have $\Delta_R(B)\subset B\otimes B$. 

\begin{thm}
  We have $\Delta(B_{\la_0}) \subset A \otimes_{<} B_{\la_0}$ and
  $\bar\Delta(B_{\la_0}) \subset B_{\la_0} \otimes_> A$. We have a
  decomposition
  $$ F = F_{2;\la_0}F_{1;\la_0},\ \on{with}\ F_{1;\la_0}\in A
  \otimes_{<} B_{\la_0}\ \on{and} \ F_{2;\la_0}\in B_{\la_0} \otimes_>
  A.
  $$ Set $\Delta_{\la_0} = \Ad(F_{1;\la_0}) \circ\Delta$, then
  $\Delta_{\la_0}$ defines a quasi-Hopf algebra structure on $A$, for
  which $B$ is a sub-quasi-Hopf algebra.
\end{thm}

{\em Sketch of proof.} The first statement is proved like Prop.\ 4.4
of \cite{ER:qH}, using the fact that $R_{\pm 2\la_0}$ are $R$-modules.
The decomposition of $F$ is proved using the same duality arguments,
e.g.\ the annihilator of $U_\hbar\N_+ \cap B_{2\la_0}$ in
$U_\hbar\N_-$ is equal to $\sum_{r\in R_{-2\la_0}}U_\hbar\N_- f[r]$.
The proof of the next statements follows \cite{ER:qH}.  \hfill \qed
\medskip

\section{Finite dimensional representations of $U_{\hbar,\omega}\G$.}

\label{fd:reps}

In \cite{Ann3}, we constructed a family $\pi_\zeta$ of $2$-dimensional
representations of $U_{\hbar,\omega}\G$ at level zero, indexed by
$\zeta$ in the infinitesimal neighborhood $\on{Spec}(\cK)$ of the
$P_i$.  We have
$$ \pi_\zeta\left( K^+(z) \right) = \pmatrix q_-(z,\zeta) & 0 \\ 0 &
q_-(q^\pa z,\zeta)^{-1} \endpmatrix,
\quad 
\pi_\zeta\left( K^-(z) \right) = \pmatrix q_+(z,\zeta) & 0 \\ 0 &
q_+(q^\pa z,\zeta)^{-1} \endpmatrix,
$$
$$ \pi_\zeta(e(z)) = \pmatrix 0 & -\hbar \sigma(z) \delta(z,\zeta) \\ 
0 & 0 \endpmatrix, \quad \pi_\zeta(f(z)) = \pmatrix 0 & 0 \\ 
\delta(z,\zeta) & 0 \endpmatrix.
$$

This family extends to a family of representations of
$U_\hbar\G^{out}$, indexed by $\zeta$ in $X - \{P_i\}$. Formulas are
$$ \pi_\zeta\left(K^+(z)\right) = \pmatrix q_-(z,\zeta) & 0 \\ 0 &
q_-(q^\pa z,\zeta)^{-1} \endpmatrix, $$
$$\pi_\zeta(e[r]) = \pmatrix 0 & -\hbar \sigma(\zeta) r(\zeta) \\ 0 &
0 \endpmatrix, \quad \pi_\zeta(f[r]) = \pmatrix 0 & 0 \\ r(\zeta) & 0
\endpmatrix.
$$ It also extends to a family of representations of
$U_\hbar\G^{out}_{\la_0}$ by the same formulas, where we fix a
preimage of $\zeta$ in $\wt X - \pi^{-1}(P_0)$. Changing the preimage
of $\zeta$ amounts to conjugating the representation by a diagonal
matrix.

Define a parenthesis order on $n$ objects as a binary tree with
extremal vertices labelled $1,\ldots,n$. To each such order, and to
$n$ points $\zeta_i$ of $X - \{P_i\}$, we associate some
$B_{\la_0}$-module.  In the case of the representation $V =
((V(\zeta_1) \otimes V(\zeta_2))\otimes (V(\zeta_3) \otimes
V(\zeta_4)))$, the space of the representation is $V = \otimes_{i=1}^n
V(\zeta_i)$ and the morphim from $B_{\la_0}$ to $\End(V)$ is
$(\otimes_{i=1}^n \pi_{\zeta_i}) \circ (\Delta\otimes\Delta) \circ
\Delta$.

In case the $\zeta_i$ are formal and $\zeta_1 << \zeta_2 << \cdots <<
\zeta_n$, the morphim $\rho_V^{(P)}$ is the restriction of a morphism
from $A$ to $\End(V)$, which is $\otimes_{i=1}^n \pi_{\zeta_i} \circ
\Ad(\Delta^{(P)}(F_{1;\la_0}))$. For example for $V = ((V(\zeta_1)
\otimes V(\zeta_2))\otimes (V(\zeta_3) \otimes V(\zeta_4)))$,
$\Delta^{(P)}(F_1)$ is equal to $\Delta^{(P)}(F_{1;\la_0}) =
F_{1;\la_0}^{(12)} F_{1;\la_0}^{(34)}(\Delta\otimes
\Delta)(F_{1;\la_0})$.

Let $(\varepsilon_1,\varepsilon_2)$ be the
canonical basis of $\CC^2$, $\xi_1,\xi_2$ its dual basis.

\begin{prop} \label{delam}
  Let $\zeta_1,\ldots,\zeta_n$ be points of $X-\{P_i\}$; let $P$ be a
  parenthesis order, and define $V$ as the $B_{\la_0}$-module
  $\otimes^{(P)}_i V(\zeta_i)$ is then a $B_{\la_0}$-module; we denote by
  $\rho_V^{(P)}$ the corresponding morphism from $B_{\la_0}$ to $\End(V)$. It
has the following properties:

1) $e[r]$ and $f[r]$ act on $V$ as $\sum_i A_i(\zeta_1,\cdots\zeta_n)
r(\zeta_i)$, $A_i$ in $\End(V)\otimes R^{\otimes n}$;

2) define the linear form $\xi$ on $V$ to be $\otimes_{i=1}^n
\xi_1^{(i)}$. Then we have $\langle \xi, \rho_V^{(P)}(f[r])v\rangle
=0$ for any $r$ of $R_{2\la_0}$ and $\langle \xi, \rho_V^{(P)}(K^+(z))
v\rangle = \prod_{i}q_+(z, \zeta_i) \langle \xi,v \rangle$, for any
$v$ in $V$. 
\end{prop}

{\em Proof.} Let us first show 1) when $\zeta_i$ are formal and
$\zeta_1<<\zeta_2<< \ldots$. In that case, $\rho_V^{(P)}(e(z))$ is conjugate
to $\Delta^{(n)}(e(z))$, which has the form 
$$\sum_i A'_i
\delta(z,\zeta_i) q_i(\zeta_i,\zeta_{i+1}, \ldots, \zeta_n),
$$ with $A'_i$ some endomorphisms of $V$ and $q_i$ in
$\CC((\zeta_{i}))\cdots((\zeta_n))$.  Therefore $\Delta^{(n)}(e[r])$
is equal some $ \sum_i A'_i r(\zeta_i) q_i(\zeta_i, \ldots, \zeta_n)$.
On the other hand, $\rho_V^{(P)}(e[r])$ is equal to the conjugation of
$\Delta^{(n)}(e[r])$ by $(\otimes_{i=1}^n
\pi_{\zeta_i})(\Delta^{(P)}(F_{1;\la_0}))$.  $\Delta^{(n)}(e[r])$
belongs to
$$\on{End}(V)((\zeta_1)) \cdots ((\zeta_n)),$$ so that
$\rho_V^{(P)}(e[r])$ has the form
\begin{equation} \label{codes}
  \sum_i B_i(\zeta_1, \ldots, \zeta_n) r(\zeta_i) , 
\end{equation}
where $B_i(\zeta_1,\ldots,\zeta_n)$ belongs to
$\End(V)((\zeta_1))\cdots ((\zeta_n))$.  $r$ being fixed,
$\rho_V^{(P)}(e[r])$ is an algebraic function in the $\zeta_i$, so the
$B_i(\zeta_1,\ldots,\zeta_n)$ are algebraic functions
and$\rho_V^{(P)}(e[r])$ has the form (\ref{codes}) for any $\zeta_i$ in
$X-\{P_i\}$. This proves 1).

Let us prove 2). Since $\Delta^{(P)}F_{1;\la_0}$ has total weight zero
(i.e. it commutes with $\Delta h[1] = \sum_i h[1]^{(i)}$),
$\rho_V^{(P)}(f[r])$ has weight $-1$, which implies the first
statement. Let us prove the second statement. We can show by induction
that
\begin{equation} \label{homer}
  \langle \xi, \Delta^{(P)}(F_{1;\la_0}) \rangle = \langle \xi,v
  \rangle,
\end{equation}
for any $v$ in $V$. For example, in the case of a representation $V =
((V(\zeta_1) \otimes V(\zeta_2))\otimes(V(\zeta_3)\otimes
V(\zeta_4)))$, we have
$$ \langle \xi, (\otimes_{i=1}^4 \pi_{\zeta_i})(F_{1;\la_0}^{(12)}
F_{1;\la_0}^{(34)}(\Delta\otimes \Delta)(F_{1;\la_0}))(v) \rangle =
\langle \xi, (\otimes_{i=1}^4 \pi_{\zeta_i})((\Delta\otimes
\Delta)(F_{1;\la_0}))(v) \rangle
$$ because $F_{1;\la_0}$ belongs to $1 + U_\hbar \N_+^{[\geq
  1]}\bar\otimes U_\hbar \N_-^{[\geq 1]}$; as
$(\Delta\otimes\Delta)(F_{1;\la_0})$ belongs to
$$1 + (U_\hbar\N_+^{[\geq 1]} \bar\otimes A^{\bar\otimes 3} )[0] +
(A\bar\otimes U_\hbar\N_+^{[\geq 1]}\bar\otimes A^{\bar\otimes 2})[0]
$$ (where $[0]$ means the zero weight component w.r.t. the adjoint
action of $\sum_i h[1]^{(i)}$), we have
$$ \langle \xi, (\Delta\otimes \Delta)(F_{1;\la_0})v\rangle =
\langle\xi,v\rangle.
$$ On the other hand, we have $\langle \xi, (\otimes_{i=1}^4
\pi_{\zeta_i})(\Delta^{(n)}(K^+(z)))(v)\rangle = \prod_i
q_+(z,\zeta_i)\langle \xi,v\rangle$. Together with (\ref{homer}), this
shows the statement for $K^+(z)$. 

\begin{remark}
  Prop.\ \ref{delam}, 2) means that the ``Drinfeld polynomial'' of
  $\otimes^{(P)}_i V(\zeta_i)$ is $\prod_{i}q_+(q^\pa z, \zeta_i)$
  (see \cite{Dr:new}). \hfill \qed \medskip
\end{remark}

Define $k_{a\to R}(z)$ as $\exp(\sum_{i,j}c_{ij} h[e^i] e^j(z))$,
where $c_{ij}$ are as in (\ref{def:ka}).

\begin{cor} \label{argive}
  There are formal series $\pi_{\al,\zeta}(z)$ such that
  $$ \langle \rho_V^{(P)}( k^+(q^{2\pa}z) k_R(q^\pa z) k_R(z)^{-1}
  k_{a\to R}(z) v) , \xi \rangle = \prod_i \pi_{\al,\zeta_i}(z)
  \langle v, \xi \rangle,
  $$ for any $v$ in $V$.
\end{cor}

\section{Twisted correlation functions} \label{tw:corr}

Let $\VV$ be a module over $U_{\hbar,\omega}\G / (K+2)$.  Let
$\psi_{\la_0}$ be a $U_\hbar\G^{out}_{\la_0}$-module map from $\VV$ to
$V$ and set $\psi_\la = \psi_{\la_0} \circ e^{\sum_a (\la_a -
  \la_a^{(0)}) h[r_a]}$.  Fix $v$ in $\VV$ and let us set
$$ f_\la(u_1,\ldots,u_n) = \langle \psi_\la [ \wt e(u_1) \cdots \wt
e(u_n) v ], \xi \rangle, 
$$ where $\xi$ is the linear form defined in Prop.\ \ref{delam}.

\begin{prop} \label{hector}
  $f_\la(u_1,\cdots,u_n)$ is a symmetric function in $(u_i)$, such
  that $$ [ \prod_{i=1}^n \prod_j \pi_{\zeta_j} (q^\pa u_i) ]
  f_\la(u_1, \cdots, u_n)$$ is regular on $\wt X^n$ except for poles
  for $u_i$ at $\pi^{-1}(P_j)$, and simple poles for $u_i$ at
  $\pi^{-1}(q^{-\pa}\zeta_j)$, and satisfy transformation properties
  (\ref{R:la}), with $\la_a^{(0)}$ replaced by $\la_a$.

  When $V$ is the trivial representation, $f_\la(u_1,\cdots,u_n)$ is
  regular on $(\wt X - \pi^{-1}(P_0))^n$.
\end{prop}

{\em Proof.} From the commutation relations of $\wt e(z)$ follows that
$f_\la(u_1,\cdots,u_n)$ is symmetric in the $u_i$. We have 
$$ \prod_j \pi_{\zeta_j}(q^\pa u_1) f_\la(u_1,\cdots,u_n) = \langle
\psi_\la[ e(q^\pa u_1) w(u_1,\cdots,u_n)] , \xi \rangle . 
$$ The fact that $ \langle \psi_\la[e[r]w] , \xi \rangle = 0$ for $r$
in $R_{2\la}$ vanishing at $\zeta_i$ implies that $$\prod_j
\pi_{\zeta_j}(q^\pa u_1) f_\la(u_1,\cdots,u_n)$$ belongs to
[(annihilator for $\langle , \rangle_\cK$ of $\{r\in R_{2\la} |
r(\zeta_i) = 0\}) \otimes \cK^{n-1}][[\hbar]]$. This annihilator is
the space of functions on $\wt X$ with simple poles at $\zeta_i$ and a
pole at the $P_i$, satisfying (\ref{R:la}). \hfill \qed \medskip

\section{Action of $T(z)$ on correlation functions} \label{F:U}

Let us set
$$
q_\mm(z,w) = (k_\mm(z), \wt e(w)), \quad \kappa(z) =
{{\exp(2\al(q^{-\pa}z,q^{-\pa}z))}\over{\exp(2\al(q^{\pa}z,q^{-\pa}z))}}
(k_a(q^{-\pa}z) , k^-(q^{-\pa}z))^{-1}.
$$

\begin{lemma}
We have
\begin{equation} \label{q:mm} 
  q_\mm(z,w) = \exp[2\hbar \sum_i ({1\over{ 1 + q^{-\pa}}} e^i)(z)
  ((T+U)e_i)_\mm(w)],
\end{equation} 
\begin{align} \label{kappa} 
  \kappa(z) = & \exp[2\al(z,z) + 2 \al(q^{-\pa}z,q^{-\pa}z) - 2
  \al(q^\pa z,z) - 2 \al(q^{-\pa}z,z)] \\ & \nonumber \exp[2\hbar \sum_i
  ({1\over{1+q^{-\pa}}} e^i)(z) ((q^\pa -1)(T+U)e_i)_{R_{(a)}}(z)]
  (k^-(z), k^-(q^{-\pa}z)) .
\end{align} 
$q_\mm(z,w)$ has the expansion
$$
q_\mm(z,w) = i_\mm(z,w) {{q^\pa z - w}\over{z-w}}, 
$$ with $i_\mm(z,w)$ in $\CC[[z,w]][z^{-1},w^{-1}][[\hbar]]^\times$.
\end{lemma}
{\em Proof.} $(k_\mm(z),\wt e(w))$ is equal to $(k_a(z),\wt
e(w))^{-1}(k^-(z),\wt e(w))$. We have already seen that $(k_a(z),\wt
e(w)) = (k_a(z), \wt f(w))^{-1}$. Then
\begin{align*}
  (k_a(z), \wt e(w)) & = (k_a(z), e(q^\pa w)) = (k_a(z),f(q^\pa
  w))^{-1} \\ & = (k_a(z),\wt f(q^\pa w) k^-(w)^{-1})^{-1} = (k_a(z),
  \wt f(w))^{-1}. 
\end{align*}
Therefore $q_\mm(z,w) = (k_\mm(z), \wt f(w))^{-1}$ and by (\ref{ccc}),
we get the statement on $q_\mm(z,w)$.  \hfill \qed \medskip

Fix $\Pi$ in $\cK[[\hbar]]$. Let $U$ be an open subset of $\CC^g$ and
define $\cF_U$ as the space of functions $f(\la_a| u_1,\cdots,u_n)$ on
$U \times (\wt X - \pi^{-1}(P_0))^{n}$, symmetric in
$(u_1,\ldots,u_n)$ and with transformation properties (\ref{R:la}),
with $(\la_a^{(0)})$ replaced by $\la_a$. For $f$ in $\cF_U$, set
\begin{align} \label{Tz}
  & (T_z^{(\Pi)} f)(\la_a | u_1, \cdots, u_n) \\ & \nonumber = \Pi(z)
  a_\la(z) \prod_{i=1}^n q_\mm(z,u_i) f(\la_a + \hbar ({1\over{1 +
      q^{-\pa}}} \omega_a / \omega)(z)
% + \sum_i b_{ai}e^i(z)
  | u_1,\cdots, u_n ) \\ & \nonumber + \Pi(q^{-\pa}z)^{-1} b'_\la(z)
  \kappa(z) \prod_{i=1}^n q_\mm(q^{-\pa}z,u_i)^{-1} f(\la_a - \hbar
  ({1\over{1 + q^{\pa}}} \omega_a / \omega)(z)
% - \sum_i b_{ai}e^i(q^{-\pa}z)
  | u_1,\cdots,u_n) \\ & \nonumber + \sum_i - {1 \over \hbar} \Pi(u_i)
  G_{2\la}(z,q^\pa u_i) q_\mm(u_i,z) \prod_{j\neq i} q_\mm(u_i,u_j) \\ 
  & \nonumber f(\la_a + \hbar ({1\over{1 + q^{-\pa}}} \omega_a /
  \omega)(u_j)
% + \sum_i b_{ai}e^i(u_j)
 | u_1,\cdots,z,\cdots,u_n) \\ & \nonumber + \sum_i {1\over \hbar}
 \Pi(q^{-\pa}u_i)^{-1} G_{2\la}(z,q^{-\pa}u_i) \kappa(u_i)
% {{\exp(2\al(q^{-\pa}u_i,q^{-\pa}u_i))} \over
%   {\exp(2\al(q^{\pa}u_i,q^{-\pa}u_i))}} (k_a(q^{-\pa}u_i),
% k^-(q^{-\pa}u_i))^{-1}
 \nonumber q_\mm(q^{-\pa}u_i,z)^{-1} \prod_{j\neq i}
 q_\mm(q^{-\pa}u_i,u_j)^{-1} \cdot \\ & \nonumber \cdot f(\la_a -
 \hbar ({1\over{1 + q^{\pa}}} \omega_a / \omega)(u_i)
%  - \sum_i b_{ai}e^i(q^{-\pa}u_i)
 | u_1,\cdots,z,\cdots,u_n)
\end{align}
where $b'_\la(z)$ is defined by (\ref{b':la}) and in the two last
sums, $q_\mm(u_i,z)$, $q_\mm(q^{-\pa}u_i,z)$, $q_\mm(u_i,u_j)$and
$q_\mm(q^{-\pa}u_i,u_j)$, $j<i$ are continued to the domains $u_i <<
z$ and $u_i << u_j$.

\begin{prop}
  Assume that $K$ acts by $-2$ on $\VV$ and $v$ is such thay $h[1]v =
  -2nv$, $\wt h[\eps] v =0$ for $\eps$ in $\mm$, and $\wt
  f[z^{1-g+k}]v = 0$ for $k\geq 0$.  We have
  $$ \langle \psi_\la[ T(z)\wt e(u_1) \cdots \wt e(u_n) v], \xi
  \rangle = T_z^{(\Pi)} (\langle \psi_\la[\wt e(u_1) \cdots \wt e(u_n)
  v], \xi \rangle),
  $$ with $T^{(\Pi)}_z$ defined by (\ref{Tz}), and $\Pi(z) = \prod_i
  \pi_{\al,\zeta_i}(z)$.
\end{prop}

{\em Proof.} We have
\begin{align*} 
  \langle & \psi_\la [T(z) \wt e(u_1) \cdots \wt e(u_n) v], \xi
  \rangle \\ & = \sum_i \langle \psi_\la [ \wt e(z) \wt e(u_1) \cdots
  [\wt f(z) , \wt e(u_i)]_{z\to R_{2\la}} \cdots \wt e(u_n) v], \xi
  \rangle \\ & + a_\la(z) \langle \psi_\la [ k^+_{tot}(z) \wt e(u_1)
  \cdots \wt e(u_n) v], \xi \rangle + b_\la(z) \langle \psi_\la [
  k^-_{tot}(z) \wt e(u_1) \cdots \wt e(u_n) v ], \xi \rangle,
\end{align*}
by the invariance of $\psi_\la$.

The sum is equal to 
\begin{align*} 
  \sum_i & -{1\over\hbar} G_{2\la}(z,q^\pa u_i) \langle \psi_\la [ \wt
  e(z) \wt e(u_1) \cdots k^+_{tot}(u_i) \cdots \wt e(u_n) v ], \xi
  \rangle \\ & + {1\over \hbar} G_{2\la}(z,q^{-\pa}u_i)
  {{\exp(2\al(q^{-\pa}u_i,q^{-\pa}u_i))} \over{\exp(2\al(q^\pa u_i,
      q^{-\pa}u_i))}} \langle \psi_\la [ \wt e(z) \wt e(u_1) \cdots
  k^-_{tot}(u_i) \cdots \wt e(u_n) v ], \xi \rangle.
\end{align*}

Then 
\begin{align*}
  & \langle \psi_\la [ k^+_{tot}(z) \wt e(u_1) \cdots \wt e(u_n) v] ,
  \xi\rangle \\ & = \Pi(z) \langle \psi_{\la_a + \sum_i b_{ai}e^i(z)}
  [ k_\mm(z) \wt e(u_1) \cdots \wt e(u_n) v ], \xi \rangle \\ & =
  \Pi(z) \prod_{i=1}^n (k_\mm(z), \wt e(u_i)) \langle \psi_{\la_a +
    \sum_i b_{ai}e^i(z)} [ \wt e(u_1) \cdots \wt e(u_n) v], \xi
  \rangle \\ & = \Pi(z) \prod_{i=1}^n q_\mm(z,u_i) \langle \psi_{\la_a
    + \sum_i b_{ai}e^i(z)} [\wt e(u_1) \cdots \wt e(u_n) v], \xi
  \rangle ;
\end{align*} 
the second equality follows from the covariance of $\psi_\la$, the
next follows from the fact that $v$ is $U_\hbar \B^{\geq
  1-g}$-invariant; in the same way
\begin{align*}
  & \langle \psi_\la [k^-_{tot}(z) \wt e(u_1) \cdots \wt e(u_n) v],
  \xi \rangle \\ & = \Pi(q^{-\pa}z)^{-1} (k_a(q^{-\pa}z),
  k_\mm(q^{-\pa}z))^{-1} \langle \psi_{\la_a - \sum_i
    b_{ai}e^i(q^{-\pa}z)} [ k_\mm(q^{-\pa}z)^{-1} \wt e(u_1) \cdots
  \wt e(u_n) v], \xi\rangle \\ & = \Pi(q^{-\pa}z)^{-1}
  (k_a(q^{-\pa}z), k^-(q^{-\pa}z))^{-1} \prod_{i=1}^n
  (k_\mm(q^{-\pa}z), \wt e(u_i))^{-1} \\ & \langle \psi_{\la_a -
    \sum_i b_{ai}e^i(q^{-\pa}z)} [ \wt e(u_1) \cdots \wt e(u_n) v],
  \xi \rangle \\ & = \Pi(q^{-\pa}z)^{-1} (k_a(q^{-\pa}z),
  k^-(q^{-\pa}z))^{-1} \prod_{i=1}^n q_\mm(q^{-\pa}z,u_i)^{-1} \\ &
  \langle \psi_{\la_a - \sum_i b_{ai}e^i(q^{-\pa}z)} [ \wt e(u_1)
  \cdots \wt e(u_n) v], \xi \rangle ,
\end{align*} 

We have
$$ -2(1\otimes q^\pa \cA)\beta(z,w) \in (R\otimes R)[[\hbar]] -2\hbar
\sum_a ({1\over{1+q^{-\pa}}} \omega_a /\omega)(z) r_a(w)
$$ so
$$ \beta(z,w) \in \hbar \sum_a ({1\over{1+q^{-\pa}}} \omega_a
/\omega)(z) r_a(w) + (R\otimes R)[[\hbar]],
$$ therefore
$$ \sum_i b_{ai} e^i(z) = \hbar ({1\over{1+q^{-\pa}}}(\omega_a
/\omega))(z)
$$ and the Proposition follows. \hfill \qed \medskip

\begin{remark} {\it Dependence on $\al$.}   
  The operators $T_z$ depend on the choice of $\al$ through their
  coefficients $\kappa(z)$ and $q_\mm(z,w)$. Operators $T_z$
  corresponding to different choices $\al$ and $\al'$ are conjugated.
  When $\Pi(z) = 1$, the conjugation is $T_z^{(\al)} = M_{\al\al'}
  T_z^{(\al')} M_{\al\al'}^{-1}$, where
  $$ (M_{\al\al'} f)(\la_a| u_1, \cdots u_n) =\prod_{i<j} \exp[ 2(\al
  - \al')(q^\pa u_i, u_j)] f(\la_a | u_1, \cdots, u_n).
  $$
\end{remark}

\section{Commuting difference operators} \label{v;chi:n}

Define $U\G_{in}^{\geq 1-g}$ as the subalgebra of $U_{\hbar,\omega}\G$
generated by $h[1]$, the $\wt h[\eps],\eps\in\mm$, and the $\wt
f[z^{1-g+k}],k\geq 0$.

Let $\chi_n$ be the character of $U\G_{in}^{\geq 1-g}$ defined by
$\chi_n(h[1]) = -2n$, $\chi_n(\wt h[\eps])= \chi_n(\wt f[z^{1-g+k}]) =
0$, $\eps$ in $\mm$, $k\geq 0$.

Define $\VV_n$ as the $U_{\hbar,\omega}\G$-module $U_{\hbar,\omega}\G
\otimes_{U\G_{in}^{\geq 1-g}} \CC_{\chi_n}$.

\begin{prop} \label{coinvts}
  The map $\iota$ from $(\VV_n^*)^{U_{\hbar}\G^{out}_{\la_0}}$ to the
  subspace $\cF$ of $S^n\cK[[\la_a -\la_a^{(0)}]][[\hbar]]$ formed of
  the formal functions near $\la_0$, which can be continued in
  variables $u_i$ to functions on $\wt X - \pi^{-1}(\{P_i\})$ with
  transformation properties (\ref{R:la}) (with $\la_a^{(0)}$ replaced
  by $\la_a$), defined by
  $$ \psi_{\la_0} \mapsto \langle \psi_\la , \wt e(u_1) \cdots \wt
  e(u_n) v \rangle, 
  $$ is an isomorphism.
\end{prop}

{\em Proof.} $(\VV_n^*)^{U_{\hbar}\G^{out}_{\la_0}}$ is isomorphic to
the space of forms $\phi$ on $U_{\hbar,\omega}\G$ such that
$\phi(x^{out}x) = \varepsilon(x^{out}) \phi(x)$, $x^{out}$ in $U_\hbar
\G_{\la_0}^{out}$ and $\phi(xx^{in}) = \phi(x) \chi_n(x^{in})$,
$x^{in}$ in $U\G_{in}^{\geq 1-g}$.

{}From Prop.\ \ref{PBW:borel} follows that the kernel of the
product map
\begin{equation} \label{borya}
  \wt \pi : U_\hbar \G^{out}_{\la_0} \otimes \CC\langle h[r_a], \wt
  e[\eps], \eps\in\cK \rangle \otimes U\G_{in}^{\geq 1-g} \to
  U_{\hbar,\omega}\G
\end{equation} 
is spanned by the $x \wt e[r_{-2\la_0}] \otimes y \otimes z -x \otimes
\wt e[r_{-2\la_0}]y \otimes z$, $r_{-2\la_0}$ in $R_{-2\la_0}$, $x$ in
$U_\hbar \G^{out}_{\la_0}$, $y$ in $\CC\langle h[r_a], \wt e[\eps],
\eps\in\cK \rangle$, $z$ in $U\G_{in}^{\geq 1-g}$. An element of $\cF$
induces a form $\bar\phi$ on $\CC\langle h[r_a], \wt e[\eps],
\eps\in\cK \rangle$, that we extend to the left side of (\ref{borya})
by the rule $\phi(x\otimes y \otimes z) = \varepsilon(x)\phi(y)
\chi_n(z)$. The properties of the elements of $\cF$ imply that $\phi$
maps $\Ker\wt\pi$ to zero. It follows that $\iota$ is surjective.

In the same way, if $\iota(\psi_{\la_0}) =0$, then the restriction of
$\iota(\psi_{\la_0})$ to $\CC\langle h[r_a], \wt e[\eps], \eps\in\cK
\rangle$ is zero, so that $\psi_{\la_0}$ is zero.  \hfill \qed
\medskip

\begin{thm} \label{comm:ops}
  For any $\Pi(z)$ in $\cK[[\hbar]]$ and $z$ in $\on{Spec}(\cK)$, the
  operators $T^{(\Pi)}_z$ defined by (\ref{Tz}) form a commuting
  family of evaluation-difference operators, acting on
  $S^n(\cK)[[\la_a - \la_a^{(0)}]][[\hbar]]$.

  When $\Pi(z) = 1$, they form a commuting family of endomorphisms of
  $\cF_U[[\hbar]]$, where $\cF_U$ is defined in sect.\ \ref{F:U}, for
  $z$ in $X - \{ P_0\}$.  Set for $\rho =
  (\rho_\la)_{\la\in\on{Spec}\CC[[\la_a - \la_a^{(0)}]]}$ a family of
  elements of $R_{ - 2\la} \cap z^{-N}\cO$,
\begin{align} \label{f:rho}
  & (\hat f[\rho]f)(\la_a | u_1,...,u_{n+1}) \\ & \nonumber =
  \sum_{i=1}^{n+1} {{1\over\hbar}} \rho_\la(u_i) \Pi(u_i) \prod_{j\neq
    i} q_\mm(u_i,u_j) f(\la_a + {{\hbar}\over{1+q^{-\pa}}} {\omega_a /
    \omega}(u_i) | u_1, ... \check i ... u_{n+1}) \\ & \nonumber -
  \sum_{i=1}^{n+1} {1\over\hbar} \rho_\la(q^{-\pa}u_i)
  \Pi(q^{-\pa}u_i)^{-1} \kappa(u_i) \prod_{j\neq i}
  q_\mm(q^{-\pa}u_i,u_j)^{-1} f(\la_a - {\hbar \over{1+q^{\pa}}}
  {\omega_a /\omega}(u_i) | u_1, ... \check i ... u_{n+1}).
\end{align}
Then the operators $T_z^{(\Pi = 1)}$ normalize the $\hat f[\rho]$,
which means that they preserve the intersection $\cap_{\rho_\la \in
  R_{ - 2\la} \cap z^{-N}\cO}\Ker\hat f[\rho]$ for any integer $N$.
\end{thm}

{\em Proof.} When $\Pi = 1$, the operators $T_z^{(\Pi)}$ can be
identified with the action of $T(z)$ on the space of invariant forms
$(\VV^*_n)^{U_\hbar \G^{out}_{\la_0}}$, by Prop.\ \ref{coinvts}.
Therefore, they preserve this space and commute with each other.

It follows that we have the cancellations of poles 
\begin{equation*} 
  \res_{z = w} [a_\la(z) q_\mm(z,w)dz] = \res_{z = w}[{1\over\hbar}
  G_{2\la}(z,q^\pa w) q_\mm(w,z) dz],
\end{equation*}
and
\begin{equation*} 
  \res_{z = w}[b'_\la(z) \kappa(z) q_\mm(q^{-\pa}z,w)^{-1} dz] +
  \res_{z = w}[{1\over\hbar} G_{2\la}(z,q^{-\pa}w)
  q_\mm(q^{-\pa}w,z)^{-1}dz] = 0.
\end{equation*}
These relations imply that when $\Pi$ is arbitrary, $T_z^{(\Pi)}$ is a
well-defined endomorphism of $S^n(\cK)[[\hbar]]$.

Set then $\Pi^+(z) = \Pi(z)$, $\Pi^-(z) = \Pi(q^{-\pa}z)^{-1}$, and 
$$ [m(\varphi(u_i))f](\la_a | u_1,...,u_n) = \varphi(u_i)f(\la_a |
u_1,...,u_n), 
$$ for any $\varphi$ in $\cK[[\hbar]]$, and
$$ T_z^{(\Pi)} = \sum_{\eps = +,-} \Pi^{\eps}(z) A^\eps_z + \sum_{\eps
  = +,-} \sum_{i=1}^n m(\Pi^\eps(u_i)) \circ C^{\eps,(i)}_z.
$$ Comparison of arguments in $(\la_a)$ in the relation $[T_z^{(\Pi =
  1)},T_w^{(\Pi = 1)}] = 0$ yields
$$
[A^\eps_z,A^{\eps'}_w] = 0, \quad [C^{\eps,(i)}_z, C^{\eps',(j)}_w] = 0 
$$
for any $\eps,\eps'$ and if $i\neq j$ and 
$$ [A^\eps_z,C^{\eps',(j)}_w] + C^{\eps',(j)}_z C^{\eps,(j)}_w = 0. 
$$
On the other hand, we have
$$ [m(\Pi^\eps(u_i)) \circ C^{\eps,(i)}_z, m(\Pi^{\eps'}(u_j)) \circ
C^{\eps',(j)}_w ] = m(\Pi^\eps(u_i)) \circ m(\Pi^{\eps'}(u_j)) \circ
 [C^{\eps,(i)}_z,C^{\eps',(j)}_w ] = 0
 $$ for any $\eps,\eps'$, $i\neq j$, and 
\begin{align*}
  & [\Pi^\eps(z) A^\eps_z, m(\Pi^{\eps'}(u_j)) \circ C^{\eps',(j)}_w]
  + m(\Pi^{\eps'}(u_j)) \circ C^{\eps',(j)}_z \circ m(\Pi^{\eps}(u_j))
  \circ C^{\eps,(j)}_w \\ & = \pi^\eps(z) m(\Pi^{\eps'}(u_j)) \circ
  ([A^\eps_z,C^{\eps',(j)}_w] + C^{\eps',(j)}_z C^{\eps,(j)}_w) \\ & =
  0. 
\end{align*}
Therefore $[T^{(\Pi)}_z,T^{(\Pi)}_w] = 0$.

The statement on $\hat f[\rho]$ follows from the fact that
$\cap_{\rho_\la \in R_{-2\la} \cap z^{-N}\cO} \Ker \hat f[\rho]$ is
equal to $(\VV_{n,N}^*)^{U_\hbar \G^{out}_{\la_0}}$, where $\VV_{n,N}$
is the $U_{\hbar,\omega}\G$-module $U_{\hbar,\omega}\G
\otimes_{U\G_{in}^{\geq -N}} \CC_{\chi_n}$, $U\G_{in}^{\geq -N}$ is
the subalgebra of $U_{\hbar,\omega}\G$ generated by the $\wt
h[\eps],\eps\in\mm$, $h[1]$ and the $\wt f[z^{k}],k\geq -N$, and
$\chi_n$ is the character of this algebra defined by $\chi_n(h[1]) =
-2n$, $\chi_n(\wt h[\eps]) = \chi_n(\wt f[z^{k}]) = 0$ for $k\geq -N$
and $\eps$ in $\mm$.  \hfill \qed \medskip

\begin{remark}
  Write $k^+(q^{2\pa}z) k_R(q^\pa z) k_R(z)^{-1} k_{a\to R}(z) =
  \exp(\sum_i h[e^i] \rho_i(z))$, with $\rho_i(z)$ in $\cK[[\hbar]]$.
  If $\Pi(z)$ has the form $\exp(\sum_i \la_i \rho_i(z))$, for some
  $\la_i$ in $\CC[[\hbar]]$, then $T_z^{(\Pi)}$ may be interpreted as
  the action of $T(z)$ on some space of intertwiners.
\end{remark}

\section{Connection with hypergeometric spaces} \label{hypergeom}

In \cite{Tar-Var}, V.\ Tarasov and A.\ Varchenko proved the following
result. Let $W$ be a representation of the Yangian $Y(\SL_2)$ and let
$\xi$ be a vector of $W$ such that $l^+_{21}(z) \xi =0$, and
$l^+_{ii}(z) \xi = \pi_i(z)\xi$, $i = 1,2$, for $\pi_i(z)$ some formal
series. 

\begin{prop} \label{TV} (see \cite{Tar-Var}). 
  We can express $(l^+_{11}(z) + l^+_{22}(z)) l^+_{12}(u_1)\cdots
  l^+_{12}(u_n) \xi$ in the form 
  $$ A(z| u_1,...,u_n) l^+_{12}(u_1)\cdots l^+_{12}(u_n) \xi +
  \sum_{i=1}^n C^{(i)}(z| u_1,...,u_n) l^+_{12}(u_1)\cdots l^+_{12}(z)
  \cdots l^+_{12}(u_n) \xi; $$ the family of operators acting on
  symmetric functions of $(u_1,...,u_n)$ defined by 
  $$ \hat T_z = A(z| u_1,...,u_n) + \sum_{i=1}^n C^{(i)}(z|
  u_1,...,u_n) \circ \ev^{(i)}_z
  $$ is commutative.
\end{prop}

In this section, we will show that the operators $\hat T_z$ are
examples of the operators $T_z^{(\Pi)}$ constructed above.

Let us consider now the case $X = \CC P^1$, $\omega = dz$. We have
$\sum_i n_i P_i = 2(\infty)$. $U_{\hbar,\omega}\G$ is then a
completion of the central extension $\wh{DY}(\SL_2)$ of the double of
the Yangian $Y(\SL_2)$ of $\SL_2$. Let $x[t^n]$, $x\in
\{e,f,h\},n\in\ZZ$ be the ``new realizations'' generators of
$DY(\SL_2)$ and $l_{ij}[n],1\leq i,j \leq 2$ and $n\in \ZZ$ its
``matrix elements'' generators.

Generators $x[t^n]$ are organized in generating series $e(z),f(z)$ and
$k^{\pm}(z)$, as above; we further split $x(z)$ as the sum $x^+(z) +
x^-(z)$, with $x^+(z) = \sum_{n\geq 0} x[t^n]z^{-n-1}$, $x^-(z) =
\sum_{n< 0} x[t^n]z^{-n-1}$.  Generating series for the $l_{ij}[n]$
are $l^{+}_{ij}(z) = \sum_{n\geq 0} l_{ij}[n]z^{-n-1}$, $l^{-}_{ij}(z)
= \sum_{n< 0} l_{ij}[n]z^{-n-1}$.

We have the relations 
$$
(z-w+\hbar) e(z) e(w) = (z-w-\hbar) e(w) e(z) , 
$$
$$ k^+(z) e(w) k^+(z)^{-1} = {{z-w+\hbar}\over{z-w}} e(w), \quad
k^-(z) e(w) k^-(z)^{-1} = {{w-z+\hbar K }\over{w-z + \hbar(K+1)}}
e(w),
$$
and
$$
l^+_{12}(z) = -\hbar k^+(z)^{-1}e^+(z), \quad 
l^-_{12}(z) = -\hbar e^-(z-\hbar K) k^-(z-\hbar)
$$
(see e.g. \cite{EF:rat}). Moreover, we have 
\begin{equation} \label{wte:yg}
  \wt e(z) = k^+(z+\hbar)^{-1} e(z+\hbar).
\end{equation} 

Define $Y^{\ge 0}$ and $Y^{<0}$ as the subalgebras of $\wh{DY}(\SL_2)$
generated the $x[t^n],n\geq 0$ (resp.\ by the $x[t^n], n<0$). Let $\wh
Y^{< 0}$ be the subalgebra generated by $K$ and $Y^{<0}$.  Then
$U_\hbar\G^{out}$ is equal to $Y^{\ge 0}$. Define $\VV$ as the Weyl
module $\wh{DY}(\SL_2) \otimes_{\wh Y^{< 0}} \CC_{-2}$, where
$\CC_{-2}$ is one-dimensional module over $\wh Y^{<0}$ where all the
generators act by zero, except for $K$, which acts by $-2$.

Let $\zeta_i$ be points of $\CC$ and $V_i(\zeta_i)$ be evaluation
modules over $Y^{\geq 0}$ associated with these points; $V_i$ is
$(2\La_i + 1)$-dimensional. Define $V$ as the tensor product (for the
usual comultiplication of $Y^{\ge 0}$) of the $V_i(\zeta_i)$. Let
$\psi$ be some $Y^{\geq 0}$-module map from $\VV$ to $V$. We will view
$V^*$ as a $Y^{\geq 0}$-module by the rule
$$
\langle a\al , v\rangle = \langle \al, S(a)v\rangle
$$ for $a$ in $Y^{\geq 0}$, $v$ in $V$ and $\al$ in $V^*$, where $S$
is the antipode of $Y^{\geq 0}$.

Let $\xi$ be a highest weight linear form as in Prop.\ \ref{delam} and
let $\Omega$ be any vector of $\VV$ annihilated by the $e^-(z)$ (for
example, $\Omega$ could be the vector $1\otimes 1$ of $\VV$). We have
$$
\langle \xi, k^+(z) v \rangle  = \pi_V(z) \langle \xi,v \rangle, 
$$
with 
$$ \pi_V(z) \pi_V(z + \hbar) = \prod_i {{\zeta_i - z + \hbar(2\La_i +
    1)}\over{\zeta_i - z}}
  $$ (see \cite{CP}).

\begin{lemma} \label{corr:corr} \label{tikva}
  Let $\wt\xi$ be any linear form of $\VV$ such that
\begin{equation} \label{archeleos}
  \langle \wt\xi, k^+(z) v\rangle = \pi (z) \langle \wt\xi,v\rangle,
  \end{equation} 
  for any $v$ in $\VV$ and some $\pi(z)$ in $\CC[[z^{-1}]]$.  Then we
  have
  \begin{equation} \label{agamemnon}
    \langle \wt\xi , e(z_1)\cdots e(z_n) \Omega \rangle =
    {1\over{(-\hbar)^n}} \prod_{i<j}{{z_j - z_i}\over{z_j - z_i
        - \hbar}} \pi(z_1) \cdots \pi(z_n) \langle \wt\xi ,
    l_{12}^+(z_1) \cdots l^+_{12}(z_n) \Omega \rangle
\end{equation} (identity in $\CC((z_1))\cdots ((z_n))$).  In particular,
  we have
  \begin{align} \label{ulysse}
    & \langle \psi(e(z_1)\cdots e(z_n) \Omega) , \xi \rangle \\ &
    \nonumber = {1\over{\hbar^n}} \prod_{i<j}{{z_j - z_i}\over{z_j -
        z_i - \hbar}} \pi_V(z_1) \cdots \pi_V(z_n)\langle
    \psi(\Omega), l^+_{12}(z_1-\hbar)\cdots l^+_{12}(z_n-\hbar) \xi
    \rangle .
 \end{align}
\end{lemma}

{\em Proof.} We proceed by induction. For $n=0$, the statement is
trivial. Assume we have proved it at step $n$ and let us try to prove
it at step $n+1$. Apply the statement of step $n$ for $\wt\xi' =
\wt\xi \circ e(z_0)$.  $\wt\xi'$ satisfies (\ref{archeleos}) with
$\pi(z)$ replaced by $\pi(z) {{z - z_0}\over{z - z_0 - \hbar}}$.
Therefore, we have
\begin{align*} 
  & \langle \wt\xi , e(z_0)\cdots e(z_n) \Omega \rangle = \langle
  \wt\xi' , e(z_1)\cdots e(z_n) \Omega \rangle \\ & =
  {1\over{(-\hbar)^n}} \pi(z_1)\cdots \pi(z_n) \prod_{0\leq i<j \leq
    n}{{z_j - z_i} \over{z_j - z_i - \hbar}} \langle \wt\xi' ,
  l_{12}^+(z_1)\cdots l_{12}^+(z_n) \Omega \rangle \\ & =
  {1\over{\hbar^n}}\pi(z_1)\cdots \pi(z_n) \prod_{0\leq i<j \leq
    n}{{z_j - z_i} \over{z_j - z_i - \hbar}} \langle \wt\xi , e(z_0)
  l_{12}^+(z_1)\cdots l_{12}^+(z_n) \Omega \rangle.
\end{align*}

Now 
\begin{align*}
  & \langle \wt\xi , e(z_0) l_{12}^+(z_1)\cdots l_{12}^+(z_n) \Omega
  \rangle \\ & = - {1\over\hbar}\langle \wt\xi , \left( k^+(z_0)
    l_{12}^+(z_0) + k^-(z_0) l_{12}^-(z_0) \right) l_{12}^+(z_1)\cdots
  l_{12}^+(z_n) \Omega \rangle \\ & = - {1\over\hbar} [ \pi(z_0)
  \langle \wt\xi ,l_{12}^+(z_0) \cdots l_{12}^+(z_n) \Omega \rangle +
  \langle \wt\xi , k^-(z_0) l_{12}^+(z_1)\cdots l_{12}^+(z_n)
  l_{12}^-(z_0) \Omega \rangle ] .
\end{align*} 
because $l_{12}^-(z_0)$ commutes with the $l_{12}^+(z_i)$. Since
$e^-(z_0) \Omega = 0$, we have $l_{12}^-(z_0) \Omega =0$, which proves
(\ref{agamemnon}) at step $n+1$.  This shows (\ref{agamemnon}).

Let us now show how (\ref{ulysse}) can be derived from
(\ref{agamemnon}). Let us set $\wt\xi(v) = \langle \xi,
\psi(v)\rangle$. Then we have (\ref{archeleos}) with $\pi(z) =
\pi_V(z)$. 
Then
\begin{align*}
  & \langle \psi(e(z_1) \cdots e(z_n) \Omega) , \xi \rangle \\ & =
  {1\over{(-\hbar)^n}} \prod_{i<j}{{z_j - z_i}\over{z_j - z_i -
      \hbar}} \pi_V(z_1) \cdots \pi_V(z_n) \langle \psi(
  l_{12}^+(z_1)\cdots l_{12}^+(z_n)\Omega) , \xi \rangle \\ & =
  {1\over{(-\hbar)^n}} \prod_{i<j}{{z_j - z_i}\over{z_j - z_i -
      \hbar}} \pi_V(z_1) \cdots \pi_V(z_n) \langle l_{12}^+(z_1)\cdots
  l_{12}^+(z_n)\psi(\Omega) , \xi \rangle \\ & = {1\over{\hbar^n}}
  \prod_{i<j}{{z_j - z_i}\over{z_j - z_i - \hbar}}\pi_V(z_1) \cdots
  \pi_V(z_n) \langle \psi(\Omega) , l_{12}^+(z_1-\hbar)\cdots
  l_{12}^+(z_n-\hbar)\xi \rangle .
\end{align*}
(the first equality by (\ref{agamemnon}); the second equality follows
from the fact that $\psi$ is a $Y^{\geq 0}$-map; the third equality
follows by definition of action on $V^*$ and because $S(l^+_{12}(z)))
= - l^+_{12}(z-\hbar)$.  \hfill \qed \medskip

\begin{cor}
We have 
$$ \langle \psi[\wt e(z_1)\cdots \wt e(z_n) \Omega], \xi \rangle =
{1\over{\hbar^n}}\langle \psi(\Omega) , l^+_{12}(z_1)\cdots
l^+_{12}(z_n) \xi \rangle.
$$
\end{cor}

{\em Proof.} We have 
\begin{align*}
  & \langle \psi[\wt e(z_1) \cdots \wt e(z_n) \Omega] , \xi \rangle \\ 
  & = \langle \psi[k^+(z_1 + \hbar)^{-1}e(z_1 + \hbar)\cdots k^+(z_n +
  \hbar)^{-1}e(z_n + \hbar) \Omega] , \xi \rangle \\ & = \prod_{i < j}
  (e(z_i), k^+(z_j)^{-1}) \prod_i \pi_V(z_i + \hbar)^{-1} \langle
  \psi[e(z_1 + \hbar)\cdots e(z_n + \hbar) \Omega] , \xi \rangle \\ &
  = {1\over{\hbar^n}} \prod_{i < j} (e(z_i), k^+(z_j)^{-1})
  \prod_{i<j} {{z_j - z_i}\over{z_j - z_i - \hbar}} \langle
  \psi(\Omega) , l^+_{12}(z_1)\cdots l^+_{12}(z_n) \xi \rangle \\ & =
  {1\over{\hbar^n}}\langle \psi(\Omega) , l^+_{12}(z_1)\cdots
  l^+_{12}(z_n) \xi \rangle ,
\end{align*}
where the first equality follows from (\ref{wte:yg}), the second from
the commutation rules, the next from Lemma \ref{tikva}.  \hfill \qed
\medskip

On the other hand, we have 
\begin{align} \label{godunov}
  & \langle \psi(T(z)\wt e(z_1)\cdots \wt e(z_n)v), \xi \rangle \\ &
  \nonumber = \langle \psi(\wt e(z_1)\cdots \wt e(z_n)T(z)v), \xi
  \rangle \on{\ (by\ centrality\ of\ }T(z)) \\ & \nonumber =
  {1\over{\hbar^n}}\langle \psi(T(z) \Omega) , l^+_{12}(z_1 + \hbar)
  \cdots l^+_{12}(z_n+\hbar) \xi \rangle
\end{align}
(by Lemma \ref{tikva} above and because $T(z)\Omega$ is killed by
$e^-(z)$).

Set $L^\pm(z) = (l^\pm_{ij}(z))_{1\leq i,j\leq 2}$, then we have $T(z)
= \on{tr} L^+(z) L^-(z - 2 \hbar) $(see \cite{RS}). But since
$l_{ij}^-(z) v = \delta_{ij}v$, we get $T(z)v = (l^+_{11}(z) +
l^+_{22}(z))v$; therefore the right side of (\ref{godunov}) is equal
to
\begin{align*}
  & {1\over{\hbar^n}} \langle \psi[(l^+_{11}(z)+ l^+_{22}(z) )\Omega],
  l^+_{12}(z_1) \cdots l^+_{12}(z_n ) \xi\rangle \\ & =
  {1\over{\hbar^n}}\langle \psi(\Omega), (l^+_{11}(z)+ l^+_{22}(z)
  l^+_{12}(z_1) \cdots l^+_{12}(z_n) \xi\rangle \\ & =
  {1\over{\hbar^n}}\hat T_z\{ \langle \psi(\Omega) ,l^+_{12}(z_1)
  \cdots l^+_{12}(z_n ) \xi\rangle \}
\end{align*}
by Prop.\ \ref{TV}.  On the other hand, $ \langle \psi[T(z)\wt
e(z_1)\cdots \wt e(z_n)\Omega], \xi \rangle$ is equal to
$$T_z^{(\Pi)}\{ \langle \psi(\Omega) ,l^+_{12}(z_1) \cdots
l^+_{12}(z_n ) \xi\rangle \} = {1\over{\hbar^n}} T_z^{(\Pi)}\{ \langle
\psi(\Omega), l^+_{12}(z_1)\cdots l^+_{12}(z_n) \xi \rangle\}
$$ by Thm.\ \ref{comm:ops}. Since any symmetric polynomial can be
realized as a correlation function $\langle \psi(\Omega),
l^+_{12}(z_1) \cdots l^+_{12}(z_n) \xi \rangle$, we have shown: 

\begin{prop}
  The operators $T_z^{(\Pi)}$ and $\hat T_z$ are equal.
\end{prop}

This fact can also be verified by direct computation.

\begin{remark} {\it Elliptic case.}
\label{elliptic}
In the elliptic case, and when there is no $z_i$, $T_z^{(\Pi = 1)}$ is
independent on $z$ and coincides with the $q$-Lam\'e operator:
$$ \hbar\theta(\hbar)(T_zf)(\la) = {{\theta(2\la -
    \hbar)}\over{\theta(2\la)}} f(\la - {\hbar\over 2}) +
{{\theta(2\la + \hbar)}\over{\theta(2\la)}} f(\la + {\hbar\over 2}). 
$$ It should be possible to obtain the $q$-Lam\'e operator for $m>1$
with other $\Pi$.

%We consider a morphism $\psi$ from $\VV$ to $V(z) = \CC^{2m+1}$, and
%set $\psi_\la(v_{top}) = \psi(\la) v[0]$. We have an equation for
%$\langle \psi_\la(e(z_1) \cdots e(z_m)v_{top}), \xi_{top}\rangle$. We
%would like to express it in terms of $\psi(\la)$. If this can be done,
%we inject in in (\ref{Tz}) and hope to get $q$-Lam\'e...

\end{remark}

\appendix

\section{Delta-function identities}

\begin{lemma} \label{rho:si}
  We have
\begin{equation} \label{meir} 
q_-(z,w)^{-1} - q_+(z,w)^{-1} = \sigma(z) \delta(q^{-\pa}z,w)
\end{equation} 
and
\begin{equation} \label{golda}
q_-(q^{-\pa}z,w) - q_+(q^{-\pa}z,w) = - \sigma(z) \delta(q^{-\pa}z,w) , 
\end{equation} 
with $\sigma$ defined by (\ref{sigma}).

$\sigma$ has also the expression 
\begin{equation} \label{rho}
 \sigma(q^\pa z) = \left[ e^{2\sum_i (U_+e_i)(z) \otimes e^i(w)}
  e^{-\phi(-\hbar,\pa_z^i \gamma)} \psi(-\hbar,\pa_z^i \gamma)
\right]_{w=z}
\end{equation}
\end{lemma}

{\em Proof.} 
From (\ref{id:q+}) follows that
$$ 
q_-(q^\pa z, w)^{-1} 
-  q_+(q^\pa z,w)^{-1} = - e^{-2\sum_i (q^{\pa}U_+ e_i)(z) \otimes
    e^i(w)}  e^{-\phi(\hbar,\pa_z^i \gamma)}
  \psi(\hbar,\pa_z^i\gamma) \delta(z,w) , 
$$
so that (\ref{meir}) follows, with $\sigma$ given by (\ref{sigma}). 

Recall that we have 
$$ q(z,w) = i(z,w) {{q^{-\pa}z-w}\over{z-q^{-\pa}w}}, 
$$ with $i(z,w)$ in $\CC[[z,w]][z^{-1},w^{-1}][[\hbar]]^\times$ such that
$i(z,w)i(w,z)=1$ (\cite{examples}, Prop.\ 3.1). We have seen that 
$$
q_-(z,w) = i_+(z,w){{q^{-\pa}z-w}\over{z-w}}, 
$$ with $i_+(z,w)$ in $\CC[[z,w]][z^{-1},w^{-1}][[\hbar]]^\times$.
Moreover, $i_+(z,w)$ satisfies
\begin{equation} \label{monory}
i_+(z,w)i_+(q^\pa z,w)={{q^{\pa}z-w}\over{z-q^{-\pa}w}}i(z,w). 
\end{equation}
On the other hand, we have 
$$
q_+(z,w) = i_+(z,w){{w-q^{-\pa}z}\over{w-z}}, 
$$
so that 
$$ q_-(q^\pa z,w)^{-1}- q_+(q^\pa z,w)^{-1} = i_+(q^\pa z,z)(q^{\pa}z-z)
\delta(z-w),
$$
and
$$
q_-(z,w) - q_+(z,w) = i_+(z,z)(q^{-\pa}z-z) \delta(z-w), 
$$ with $\delta(z-w) = \sum_{i\in\ZZ} z^i w^{-i-1}$.  From $i(z,z)=1$
and (\ref{monory}) follows that the prefactors of $\delta(z-w)$ in
both equations are opposite to each other. (\ref{golda}) follows.

On the other hand, we have 
$$ q_+(z,w) = q^{2\sum_i (U_+e_i)(z) \otimes e^i(w)}
e^{-\phi(-\hbar,\pa_z^i\gamma)} (1 + G^{(21)}(z,w) \psi(-\hbar,\pa_z^i
\gamma)), 
$$
so that 
$$
q_-(z,w) - q_+(z,w) = \rho(q^\pa z)\delta(z,w). 
$$
with $\rho$ given by 
$$
 \rho(q^\pa z) = \left[ - e^{2\sum_i (U_+e_i)(z) \otimes e^i(w)}
  e^{-\phi(-\hbar,\pa_z^i \gamma)} \psi(-\hbar,\pa_z^i \gamma)
\right]_{w=z} ; 
$$ since $\rho(q^\pa z)$ is equal to $-\sigma(q^\pa z)$, we get
expression (\ref{rho}) for $\sigma$.  \hfill \qed \medskip

\begin{lemma} \label{al:beta}
We have 
\begin{equation} \label{B}
  q_+(z,w)^{-1} G(w,z) + q_-(z,w)^{-1} G(z,w) 
= \al(z) \delta(q^{-\pa}z,w), 
\end{equation}
\begin{equation} \label{A'''}
  q_+(z,w) G(w,q^{-\pa}z) + q_-(z,w) G(q^{-\pa}z,w) = \beta(q^{\pa}z)
  \delta(z,w),
\end{equation}
with $\al$ and $\beta$ defined by (\ref{alpha}) and (\ref{beta}).
$\beta$ has also the expression
\begin{equation} \label{alt:beta} 
  \beta(q^{\pa} z) = \left[ \pa_{\hbar}[ e^{-\phi(-\hbar,\pa_z^i
      \gamma)} \psi(-\hbar,\pa_z^i\gamma)] e^{2\sum_i
      (U_+e_i)(z)\otimes e^i(w)} \right]_{w=z}.
\end{equation}
\end{lemma}

{\em Proof.}
Let us prove (\ref{B}). Applying $q^\pa \otimes 1$ to this
equation, we write it as 
$$ q_+(q^\pa z,w)^{-1} G(w,q^\pa z) + q_-(q^\pa z,w)^{-1} G(q^\pa z,w)
= \al(q^\pa z)\delta(z,w).
$$

We have
$$
q_+(q^\pa z,w)^{-1} = e^{ - q^\pa U_+ e_i(z) \otimes e^i(w)}
e^{\sum_i {{1-q^\pa}\over{\pa}}e_i(z) \otimes e^i(w)}. 
$$ 
From sect.\ \ref{kernels} follows that
\begin{align} \label{B''}
e^{\sum_i {{1-q^\pa}\over{\pa}}e_i(z) \otimes e^i(w)}
& = e^{-\phi(-\hbar,(-\pa_z)^i \gamma)}
  (1-G^{(21)}\psi(-\hbar,(-\pa_z)^i\gamma))
\\ & \nonumber 
= e^{-\phi(\hbar,\pa_z^i \gamma)}
  (1+G^{(21)}\psi(\hbar,\pa_z^i\gamma)) .
\end{align}

$q_+(q^\pa z,w)^{-1} G^{(21)}(q^\pa z,w)$ is equal to
\begin{align*} 
& q^{ - 2 q^\pa (T_+ + U_+)e_i(z) \otimes e^i(w)} G^{(21)}(q^\pa z,w)
  \\ & = 
e^{ ({{1-q^\pa}\over{\pa}} - 2 q^\pa U_+)e_i(z) \otimes e^i(w)}
  G^{(21)}(q^\pa z,w) 
\\ & 
= - e^{- 2\sum_i q^\pa U_+ e_i(z) \otimes e^i(w)}
\pa_{\hbar}(e^{\sum_i {{1-q^\pa}\over{\pa}}e_i(z) \otimes e^i(w)})
\\ & 
= - e^{-2\sum_i q^\pa U_+ e_i(z) \otimes e^i(w)}
\pa_{\hbar}[
e^{-\phi(\hbar,\pa_z^i\gamma)}(1+G^{(21)}\psi(\hbar,\pa_z^i\gamma))]
\end{align*}
therefore 
\begin{align*}
  & q_+(q^\pa z,w)^{-1} G^{(21)}(q^\pa z,w) -
  q_-(q^{-\pa}z,w)^{-1}G(q^\pa z,w) \\ & = - e^{-2 \sum_i q^\pa U_+
    e_i(z) \otimes e^i(w)}
  \pa_{\hbar}[e^{-\phi(\hbar,\pa_z^i\gamma)}\psi(\hbar,\pa_z^i\gamma))]
  \delta(z,w).
\end{align*}
Therefore 
$$ \al(q^\pa z) = \left[ - e^{-\sum_i q^\pa U_+ e_i(z) \otimes e^i(w)}
\pa_{\hbar} \{e^{-\phi(\hbar,\pa_z^i\gamma)}\psi(\hbar,\pa_z^i\gamma))\}
\right]_{z = w}. 
$$

Let us now prove (\ref{A'''}).  From sect.\ \ref{kernels} follows that
\begin{equation} \label{A'}
e^{\sum_{i} ({{1-q^{-\pa}}\over{\pa}} e_i)(z) \otimes e^i(w)}  
=   e^{-\phi(\hbar,(-\pa_z)^i \gamma)}
  (1-G^{(21)}\psi(\hbar,(-\pa_z)^i\gamma)).
\end{equation}
Differentiating (\ref{A'}) with respect to $\hbar$, we find 
$$
G^{(21)}(q^{-\pa}z,w) e^{\sum_{i} ({{1-q^{-\pa}}\over{\pa}} e_i)(z)
\otimes e^i(w)} = \pa_{\hbar}[ e^{-\phi(\hbar,(-\pa_z)^i \gamma)}
(1-G^{(21)}\psi(\hbar,(-\pa_z)^i\gamma))].
$$
Therefore, we have also 
$$
G(q^{-\pa}z,w) (e^{\sum_{i} ({{1-q^{-\pa}}\over{\pa}} e_i)(z)
\otimes e^i(w)})_{w<< z} 
= \pa_{\hbar}[ e^{-\phi(\hbar,(-\pa_z)^i \gamma)}
(-1-G\psi(\hbar,(-\pa_z)^i\gamma))] , 
$$
so that 
\begin{align*}
& G^{(21)}(q^{-\pa}z,w) e^{\sum_{i} ({{1-q^{-\pa}}\over{\pa}} e_i)(z)
\otimes e^i(w)} 
- 
[G^{(21)}(q^{-\pa}z,w) e^{\sum_{i} ({{1-q^{-\pa}}\over{\pa}} e_i)(z)
\otimes e^i(w)} ]_{w << z} 
\\ & = - \pa_{\hbar}[ e^{-\phi(\hbar,(-\pa_z)^i \gamma)}
\psi(\hbar,(-\pa_z)^i\gamma)] \delta(z,w), 
\end{align*}
and
\begin{align*}
& G^{(21)}(q^{-\pa}z,w) q_+(z,w)+ G(q^{-\pa}z,w) q_-(z,w)
\\ & = - \pa_{\hbar}[ e^{-\phi(\hbar,(-\pa_z)^i \gamma)}
\psi(\hbar,(-\pa_z)^i\gamma)] e^{2\sum_i (U_+e_i)(z)\otimes e^i(w)}
\delta(z,w),  
\end{align*}
that is (\ref{A'''}) with $\beta$ given by (\ref{alt:beta}). Identity
(\ref{rho}) allows then to write $\beta$ in the form (\ref{beta}).
\hfill \qed \medskip

\begin{lemma} \label{ids:G:q}
  We have
  $$ G_{-2\la}(q^\pa w,z) q_-(q^{2\pa}w,q^\pa z)^{-1} +
  G_{2\la}(z,q^\pa w) q_+(q^{2\pa}w,q^\pa z)^{-1} = A_\la(z)
  \delta(z,w),
  $$ and
  $$ G_{-2\la}(q^{-\pa} w,z) q_-(q^\pa w, q^\pa z) +
  G_{2\la}(z,q^{-\pa} w) q_+(q^\pa w, q^\pa z) = B_\la(z) \delta(z,w), 
  $$ where $A_\la(z)$ and $B_\la(z)$ are defined by (\ref{A:la}) and
  (\ref{B:la}).
\end{lemma}

{\em Proof.}
Using Lemmas \ref{al:beta} and \ref{rho:si}, we find
$$ A_\la(z) = \al(q^{2\la}z) + [G_{-2\la}(q^\pa w,z) -
G(q^{2\pa}w,q^\pa z)]_{| w = z} \sigma(q^{2\pa}z). 
$$
Then $G(q^\pa w,z) - G(q^{2\pa}w,q^\pa z)$ is equal to 
\begin{equation} \label{tensor}
\sum_i e_i \otimes q^\pa e^i - \sum_i q^\pa e_i \otimes q^{2\pa} e^i ; 
\end{equation}
the pairing of $R$ with the first component of this tensor gives zero,
so that it belongs to $(R\otimes R)[[\hbar]]$ and its pairing with
$\la$ in $\La$ gives $-q^{2\pa}(q^{-\pa}\la)_R$; therefore
(\ref{tensor}) is equal to
$$ - \sum_i e^i \otimes q^{2\pa}(q^{-\pa}e_i)_R.
$$ so that $A_{\la}(z)$ is given by (\ref{A:la}).

In the same way, we find
$$ B_\la(z) = \beta(q^{2\pa}z) - \sigma(q^{2\pa}z)
[G_{-2\la}(q^{-\pa}w,z)- G(w,q^\pa z) ]_{|w = z}. 
$$
Since 
$$
G(q^{-\pa}w,z)- G(w,q^\pa z) 
$$ is equal to $ - \sum_i e^i(z) ((q^{-\pa}e_i)_R)(z)$, it follows
that $B_\la(z)$ is given by (\ref{B:la}).  \hfill \qed 

\end{document}